\documentclass{article}
\title{On Cores in Yetter-Drinfel'd Hopf Algebras}
\author{Yevgenia Kashina, Yorck Sommerh\"auser}
\date{}

\usepackage{amsmath}
\usepackage{amsthm}
\usepackage{amssymb}
\usepackage{latexsym}
\usepackage{array}
\usepackage{mathtools}

\makeatletter
\renewcommand{\subsection}{\@startsection{subsection}{2}{0em}%
{\baselineskip}{-0em}{\bfseries\normalsize}}
\makeatother

\newcounter{num}
\newenvironment{pflist}{\begin{list}{(\arabic{num})}{\usecounter{num} \leftmargin0cm \itemindent5pt}}{\end{list}}

\newenvironment{thmlist}{\begin{list}{\arabic{num}.}{\usecounter{num} \leftmargin0pt \itemindent5pt \topsep-5pt \partopsep0pt \labelwidth0pt}}{\end{list}}

\makeatletter
\def\enumfix{%
\if@inlabel
 \noindent \par\nobreak\vskip-\topsep\hrule\@height\z@
\fi}
\let\oldenumerate\enumerate
\def\enumerate{\enumfix\oldenumerate}
\let\oldthmlist\thmlist
\def\thmlist{\enumfix\oldthmlist}
\makeatother

\theoremstyle{definition}
\newtheorem*{defn}{Definition}
\newtheorem*{prop}{Proposition}
\newtheorem*{lem}{Lemma}
\newtheorem*{cor}{Corollary}
\newtheorem*{conj}{Conjecture}
\newtheorem*{pf}{Proof}

\newcommand{\dm}{\operatorname{dim}}

\newcommand{\Aut}{\operatorname{Aut}}
\newcommand{\Span}{\operatorname{Span}}
\newcommand{\Hom}{\operatorname{Hom}}
\newcommand{\id}{\operatorname{id}}

\newcommand\1{{(1)}}
\newcommand\2{{(2)}}
\newcommand\3{{(3)}}

\newcommand\ot{\otimes}

\newcommand\da{\Delta_{A}}

\newcommand\Dh{\Delta_{H}}
\newcommand\ea{\varepsilon_{A}}

\newcommand\eh{\varepsilon_{H}}

\newcommand\sa{S_{A}}

\newcommand\sh{S_{H}}

\newcommand\uh{1_H}
\newcommand\A{1_A}
\newcommand\ua{1}
\newcommand\G{1}
\newcommand\K{1}

\newcommand\Z{{\mathbb Z}}

\newcommand\deq{\vcentcolon =}
\newcommand\coln{\colon}

\begin{document}

\maketitle

\begin{abstract}
\hspace{-17pt} 
By constructing explicit examples, we show that the core of a group-like element in a cocommutative cosemisimple Yetter-Drinfel'd Hopf algebra over the group ring of a finite abe\-lian group is not always completely trivial.
\end{abstract}

\section*{Introduction} \label{Sec:Introd}
\addcontentsline{toc}{section}{Introduction}
Even if the base field is algebraically closed of characteristic zero, there is at present no complete structure theory for finite-dimensional cocommutative cosemisimple Yetter-Drinfel'd Hopf algebras over the group ring of a finite abelian group. Such a structure theory is only available if the finite abelian group has prime order (cf.~\cite{SoYp}). In the prime order case, it turns out that the so-called core of a group-like element is always completely trivial in the sense that both the action and the coaction on the core are trivial. The notion of the core itself can be defined in the case of a general finite abelian group. This raises the question whether also in this case the core of a group-like element is always completely trivial. The purpose of this article is to construct examples that show that this is not the case. 

The examples that we construct are eight-dimensional Yetter-Drinfel'd Hopf algebras over the group ring of the elementary abelian group~$\Z_2 \times \Z_2$. They have a basis consisting of group-like elements and are therefore cocommutative and cosemisimple. After a general discussion of Yetter-Drinfel'd Hopf algebras in Section~\ref{Sec:Prelim}, where we also discuss the notions of triviality and complete triviality, we construct in Section~\ref{Sec:FirstExamp} a commutative example of such a Yetter-Drinfel'd Hopf algebra. The example depends on a parameter~$\zeta$, which is a not necessarily primitive fourth root of unity. In Section~\ref{Sec:SecExamp}, we construct a second such example, which is very similar, but noncommutative. Also the second example depends on a not necessarily primitive fourth root of unity~$\zeta$. 

In Section~\ref{Sec:Cores}, we first review parts of the theory developed in~\cite{SoTriv} in a dualized form for cocommutative instead of commutative Yetter-Drinfel'd Hopf algebras. In particular, we review the notion of the core of a group-like element. We then discuss the cores that appear in the examples constructed before and in particular observe that some are not completely trivial, because the finite abelian group acts and coacts on them nontrivially. However, these cores are trivial in the sense that they are ordinary Hopf algebras, while the Yetter-Drinfel'd Hopf algebras that we constructed are nontrivial; i.e., they are not ordinary Hopf algebras. The article concludes with the conjecture that this is a general phenomenon: We conjecture that the core of a group-like element is always trivial, at least under the hypotheses that we have been making throughout the article. 

The expert reader will realize immediately that the examples that we construct are crossed products over a group with two elements, and wonder why this fact is not mentioned in the text. The reason is that much more is true: Not only the algebra structure arises from a specific construction, but also the coalgebra structure, and it is possible to describe the necessary compatibility conditions between the data involved. We plan to address these topics in the forthcoming article~\cite{KaSo}. However, as a consequence of the fact that the cores in our examples are not completely trivial, our examples are not diagonal realizations in the sense of~\cite{AndNa}, Sec.~3, a construction that can also yield crossed products over a group with two elements.

Concerning conventions, the precise assumptions that are used throughout each section are listed in its first paragraph. In particular, this holds for the assumptions on the base field, which is always denoted by~$K$. It is arbitrary in Section~\ref{Sec:Prelim}, but in Section~\ref{Sec:FirstExamp}, we assume that it contains a primitive eighth root of unity, which forces that its characteristic is different from~$2$. In Section~\ref{Sec:SecExamp}, we assume that it contains a primitive fourth root of unity. Finally, in Section~\ref{Sec:Cores}, we assume that it is algebraically closed of characteristic zero. All vector spaces are defined over~$K$, and all unadorned tensor products are taken over~$K$; more general tensor products appear only once in Paragraph~\ref{AlgStruc}. The dual of a vector space~$V$ is denoted by $V^* \deq \Hom_K(V,K)$, and the transpose of a linear map~$f$, i.e., the induced map between the dual spaces, is denoted by~$f^*$. The set of matrices with~$m$ rows, $n$~columns, and entries in~$K$ is denoted by $M(m \times n, K)$, and the $n \times n$-identity matrix is denoted by $E_n \in M(n \times n, K)$.

\enlargethispage{8pt}
\begin{samepage}
All algebras are assumed to have a unit element, and algebra homomorphisms are assumed to preserve unit elements. Unless stated otherwise, a module is a left module. The multiplicative group of invertible elements in the base field~$K$ is denoted by~$K^\times \deq K\setminus\{0\}$.  In this article, a character is a one-dimensional character, i.e., a group homomorphism to the multiplicative group~$K^\times$ in the case of a group character, or an algebra homomorphism to the base field~$K$ in the case of the character of an algebra. The character group of an abelian group~$G$ is the group of all such characters, i.e., the group~$\hat{G} \deq \Hom(G,K^\times)$. The symbol~$\perp$ will not be used for vector spaces, but only for abelian groups. It~has two meanings: For a subset of the group~$G$, it denotes the set of characters that take the value~$\K$ on all elements of the set, and for a subset of the character group~$\hat{G}$, it denotes the set of all group elements on which all of the given characters take the value~$\K$. 
\end{samepage}

The cardinality of a set~$X$ will be denoted by~$|X|$. The symbol~$\subset$ denotes non-strict inclusion, so that we have $X \subset X$ for every set~$X$. Also, we use the so-called Kronecker symbol~$\delta_{ij}$, which is equal to~$1$ if $i=j$ and equal to~$0$ otherwise. With respect to enumeration, we use the convention that propositions, definitions, and similar items are referenced by the paragraph in which they occur; i.e., a reference to Proposition~1.1 refers to the unique proposition in Paragraph~1.1. 

The material discussed here was first presented at the AMS Spring Southeastern Sectional Meeting in Auburn in March~2019. The travel of the second author to this conference as well as his work on this article were supported by NSERC grant RGPIN-2017-06543. The work of the first author on this article was supported by a Faculty Summer Research Grant from the College of Science and Health at DePaul University.

\section{Preliminaries} \label{Sec:Prelim}
\subsection[Yetter-Drinfel'd modules]{} \label{YetMod}
In this section, we review some basic facts about Yetter-Drinfel'd Hopf algebras that will be needed in the sequel. Our base field~$K$ can be arbitrary. We consider a Hopf algebra~$H$ over~$K$ with coproduct~$\Dh$, counit~$\eh$, and bijective antipode~$\sh$ (cf.~\cite{M}, Def.~1.5.1, p.~7). A left-left Yetter-Drinfel'd module is a left $H$-module~$V$ that is simultaneously a left $H$-comodule in such a way that these two structures are compatible in the sense that
\[\delta(h.v) = h_{\1} v^{\1} \sh(h_{\3}) \ot h_{\2}.v^{\2}\]
for all $h \in H$ and all $v \in V$. Here we have used Heyneman-Sweedler sigma notation in the form 
\[\Dh(h) = h_{\1} \ot h_{\2} \in H \ot H \qquad \qquad
\delta(v) = v^{\1} \ot v^{\2} \in H \ot V\]
for the coproduct~$\Dh$ and the coaction~$\delta$, respectively, and a dot to denote the module action. Although there are also left-right, right-left, and right-right Yetter-Drinfel'd modules, by a Yetter-Drinfel'd module we always mean a left-left Yetter-Drinfel'd module unless specified otherwise. When~$H$ is commutative and cocommutative, the Yetter-Drinfel'd condition stated above reduces to
\[\delta(h.v) = v^{\1} \ot h.v^{\2}\]
which coincides with the compatibility condition for left-left dimodules (cf.~\cite{Long}, Def.~3.1, p.~575).

The tensor product of two Yetter-Drinfel'd modules, together with the usual diagonal module structure and the codiagonal comodule structure, is again a Yetter-Drinfel'd module. If~$V$ and~$W$ are two Yetter-Drinfel'd modules, the two possible ways of forming their tensor product are related via the quasisymmetry
\[\sigma_{V,W} \coln V \ot W \rightarrow W \ot V,~v \ot w \mapsto (v^{\1}.w \ot v^{\2})\]
This map is both $H$-linear and colinear. Since we are assuming that the antipode of~$H$ is bijective, the quasisymmetry is also bijective, and its inverse is given by the formula 
$\sigma_{V,W}^{-1}(w \ot v) = v^{\2} \ot \sh^{-1}(v^{\1}).w$. We say that the action on~$V$ is trivial if we have $h.v= \eh(h) v$ for all $h \in H$ and all~$v \in V$, i.e., if every element of~$V$ is invariant, and we say that the coaction on~$V$ is trivial if we have $\delta(v) = \uh \ot v$ for all~$v \in V$, i.e., if every element of~$V$ is coinvariant. Any vector space~$V$, and in particular the base field~$K$, becomes a Yetter-Drinfel'd module if endowed with the trivial action and the trivial coaction. We will always view the base field as a Yetter-Drinfel'd module in this way. Note that the quasisymmetry coincides with the flip, i.e., we have
\[\sigma_{V,W}(v \ot w) = w \ot v \]
for all~$v \in V$ and all~$w\in W$, if the coaction on~$V$ is trivial or the action on~$W$ is trivial.

\subsection[Yetter-Drinfel'd Hopf algebras]{} \label{YetDrinfHopf}
What we have said in the preceding paragraph can be extended to the statement that Yetter-Drinfel'd modules form a quasisymmetric monoidal category. While algebras and coalgebras can be defined in any monoidal category, the quasisymmetry is necessary to define bialgebras and Hopf algebras. By Yetter-Drinfel'd algebras, Yetter-Drinfel'd coalgebras, Yetter-Drinfel'd bialgebras, and Yetter-Drinfel'd Hopf algebras, we mean the respective algebraic object in the category of Yetter-Drinfel'd modules. A Yetter-Drinfel'd Hopf algebra~$A$ therefore comes endowed with five structure elements: 
\begin{enumerate}
\item 
a multiplication $\mu_A \coln A \ot A \to A$, denoted by $\mu_A(a \ot a') = aa'$,

\item
a unit map $\eta_A \coln K \to A$, used to define $\A \deq \eta_A(\K)$,

\item 
a coproduct~$\da \coln A \to A \ot A$, for which we use, as for the Hopf algebra~$H$, the notation
$\da(a) = a_{\1} \ot a_{\2}$,

\item
a counit $\ea \coln A \to K$,

\item 
and an antipode~$\sa \coln A \to A$.
\end{enumerate}
All of these five structure elements must be morphisms in the category; i.e., they need to be $H$-linear and colinear. When endowed with these structure elements, a Yetter-Drinfel'd Hopf algebra is usually not an ordinary Hopf algebra, because in the requirement that the coproduct be an algebra homomorphism, the algebra structure on~$A \ot A$ is given by
\[\mu_{A \ot A} \deq (\mu_{A} \ot \mu_{A}) \circ (\id_A \ot \, \sigma_{A,A} \ot \id_{A})\]
and therefore involves the quasisymmetry~$\sigma_{A,A}$ at the place where the canonical tensor product would involve the flip. We denote this algebra structure by~$A \hat{\ot} A$ to distinguish it from the canonical tensor product algebra structure. 

On the other hand, this means that a Yetter-Drinfel'd Hopf algebra is an ordinary Hopf algebra, endowed with additional structure, if its quasisymmetry~$\sigma_{A,A}$ coincides with the flip, i.e., if we have
\[\sigma_{A,A}(a \ot a') = a' \ot a\]
for all~$a, a' \in A$. As pointed out above, this will happen if the action or the coaction of~$H$ on~$A$ are trivial, and in particular if both the action and the coaction are trivial. We formalize this in the following definition:
\begin{defn} 
\leavevmode 
\begin{thmlist}
\item 
A Yetter-Drinfel'd Hopf algebra~$A$ is called trivial if we have 
\[\sigma_{A,A}(a \ot a') = a' \ot a\]
for all $a,a' \in A$.

\item
A Yetter-Drinfel'd Hopf algebra~$A$ is called completely trivial if both the action and the coaction of~$H$ on~$A$ are trivial.
\end{thmlist}
\end{defn}

The first definition should be compared with the ones made in~\cite{SoDiss}, Def.~1.1, p.~9, \cite{SoYp}, Def.~1.1, p.~8, and \cite{SoTriv}, Def.~1.3, p.~481. Our statement above now says that a trivial Yetter-Drinfel'd Hopf algebra is an ordinary Hopf algebra, endowed with some additional structure. As observed by P.~Schauenburg (cf.~\cite{Schau}, Cor.~2, p.~262), the converse of this assertion is also true; i.e., a Yetter-Drinfel'd Hopf algebra that is, with respect to the same structure elements, an ordinary Hopf algebra must necessarily be trivial. Short proofs of this fact in our setting are given in~\cite{SoDiss}, Prop.~1.1, p.~9 or~\cite{SoYp}, Prop.~1.1, p.~8. On the other hand, there is no essential difference between a completely trivial Yetter-Drinfel'd Hopf algebra and an ordinary Hopf algebra. 

We have said above that all structure elements of a Yetter-Drinfel'd Hopf algebra, in particular the antipode, are required to be $H$-linear and colinear. While this is correct, it turns out that in the case of the antipode, this requirement is automatically satisfied:
\begin{lem} 
Let $A$ be a Yetter-Drinfel'd bialgebra. Assume that $\id_A$ has a convo\-lution-inverse~$\sa \coln A \to A$, i.e., a $K$-linear map that satisfies
\[\sa(a_\1) a_\2 = a_\1 \sa(a_\2) = \ea(a) \A \]
Then $\sa$ is both $H$-linear and $H$-colinear.
\end{lem}
\begin{pf} 
We have
\begin{align*}
\sa (h.a) &= \sa (h_\1.a_\1) \left(h_\2. \! \left( a_\2 \sa (a_\3) \right) \right)\\
&=\sa (h_\1.a_\1) \left( h_\2. a_\2\right)\left(h_\3.\sa (a_\3) \right) \\
&=
\ea (h_\1.a_\1) \left(h_\2.\sa (a_\2)\right) = h.\sa(a)
\end{align*}
and
\begin{align*}
a^\1 \ot \sa(&a^\2) = a_\2{}^\1 \ot \ea (a_\1) \sa(a_\2{}^\2)  \\
&= \left(\ea (a_\1) \A \right)^\1 a_\2{}^\1 \ot 
\left(\ea (a_\1) \A \right)^\2 \sa(a_\2{}^\2) \\
&= \left(\sa (a_\1) a_\2 \right)^\1 a_\3{}^\1 \ot 
\left(\sa (a_\1) a_\2 \right)^\2 \sa(a_\3{}^\2) \\
&= \sa(a_\1)^\1 a_\2{}^\1 a_\3{}^\1 \ot \sa (a_\1)^\2 a_\2{}^\2 \sa(a_\3{}^\2) \\
&= \sa(a_\1)^\1 a_\2{}^\1 \ot \sa (a_\1)^\2 a_\2{}^\2{}_\1 \sa(a_\2{}^\2{}_\2)\\
&= \sa(a_\1)^\1 a_\2{}^\1 \ot \sa (a_\1)^\2 \ea (a_\2{}^\2) 
= \sa(a)^\1 \ot \sa (a)^\2
\end{align*}
Thus $\sa$ is $H$-linear and colinear.
\qed
\end{pf}

\subsection[Quasitriangular Hopf algebras]{} \label{QuasitriHopf}
A way of constructing Yetter-Drinfel'd modules and Yetter-Drinfel'd Hopf algebras uses quasitriangular Hopf algebras. Recall that a Hopf algebra~$H$ with bijective antipode is called quasitriangular if it possesses a so-called universal R-matrix, which is an invertible element 
$R = \sum_{i=1}^m a_i \ot b_i \in H \ot H$ that satisfies
\begin{enumerate}
\item 
$\displaystyle
\sum_{i=1}^m h_\2 a_i \ot h_\1 b_i = \sum_{i=1}^m a_i h_\1 \ot b_i h_\2$ for all $h \in H$

\item 
$\displaystyle
(\Dh \ot \id)(R) = \sum_{i,j=1}^m a_i \ot a_j \ot b_i b_j$

\item 
$\displaystyle
(\id \ot \Dh)(R) = \sum_{i,j=1}^m a_i a_j \ot b_j \ot b_i$
\end{enumerate}
(cf.~\cite{M}, Def.~10.1.5, p.~180). With the help of an R-matrix, any left $H$-module~$V$ can be made into a Yetter-Drinfel'd module by introducing the coaction
\[\delta \coln V \to H \ot V,~v \mapsto \sum_{i=1}^m b_i \ot a_i.v\]
(cf.~\cite{M}, Examp.~10.6.14, p.~213). It follows from the fundamental properties of $R$-matrices (cf.~\cite{M}, Prop.~10.1.8, p.~180f) that $H$-invariant elements are also coinvariant with respect to this coaction. In particular, if the $H$-module structure on~$V$ is trivial, then the resulting comodule structure will also be trivial. 

If we turn a second $H$-module~$W$ into a Yetter-Drinfel'd module in this way, an $H$-linear map from~$V$ to~$W$ is also $H$-colinear with respect to these comodule structures. Furthermore, the formula for the quasisymmetry recorded in Paragraph~\ref{YetMod} becomes
\[\sigma_{V,W}(v \ot w) = \sum_{i=1}^m b_i.w \ot a_i.v\]
which is the standard quasisymmetry on $H$-modules induced by an R-matrix (cf.~\cite{M}, Thm.~10.4.2, p.~199). Moreover, if we consider $V \ot W$ as an $H$-module and then endow it with the coaction just described, we obtain exactly the codiagonal comodule structure that we have used in Paragraph~\ref{YetMod}. Stated in categorical terms, we have constructed a strict quasisymmetric monoidal functor from the category of $H$-modules to the category of Yetter-Drinfel'd modules (cf.~\cite{Kas}, Def.~XIII.3.6, p.~327). Therefore, a Hopf algebra in the category of \mbox{$H$-modules} will be a Yetter-Drinfel'd Hopf algebra when endowed with the above comodule structure.

\subsection[Bicharacters]{} \label{Bichar}
One source of R-matrices are bicharacters on finite abelian groups: 
\begin{defn}
If~$G$ is a finite abelian group, written multiplicatively, a function
\[\theta \coln G \times G \to K^\times,~(g,g') \mapsto \theta(g,g')\]
is called a bicharacter if we have
\[\theta(g,g'g'') = \theta(g,g') \theta(g,g'')\qquad \text{and} \qquad 
\theta(gg', g'') = \theta(g,g'') \theta(g',g'')\]
for all $g, g', g'' \in G$. The bicharacter is called symmetric if $\theta(g,g') = \theta(g',g)$. It is called nondegenerate if the only element~$g \in G$ that satisfies $\theta(g,g') = \K$ for all $g' \in G$ is the unit element $g=\G$. 
\end{defn}

The property of nondegeneracy in fact implies the corresponding condition on the other side: If~$\theta$ is nondegenerate, the only element~$g' \in G$ that satisfies $\theta(g,g') = \K$ for all $g \in G$  is the unit element $g'=\G$. 

The two equations in the following lemma are called orthogonality relations:
\begin{lem}
If~$\theta$ is a nondegenerate bicharacter, we have
\[\sum_{g \in G} \theta(g,g') = |G| \delta_{g',1} \qquad \text{and} \qquad 
\sum_{g' \in G} \theta(g,g') = |G| \delta_{g,1}\]
\end{lem}
\begin{pf}
The first assertion is obvious when $g'=\G$. If $g' \neq \G$, nondegeneracy implies that there is $h \in G$ such that $\theta(h,g') \neq \K$. We then have
\begin{align*}
\sum_{g \in G} \theta(g,g') = \sum_{g \in G} \theta(gh,g') = \sum_{g \in G} \theta(g,g') \theta(h,g')
= \theta(h,g') \sum_{g \in G} \theta(g,g')
\end{align*}
which implies $\sum_{g \in G} \theta(g,g') = 0$. The second assertion follows by applying the first assertion to the bicharacter that arises from~$\theta$ by interchanging the arguments. 
\qed
\end{pf}

\enlargethispage{7pt}
\begin{samepage}
Bicharacters lead to R-matrices in the following way:
\begin{prop}
Suppose that~$G$ is a finite abelian group whose order~$|G|$ is not divisible by the characteristic of~$K$. If~$\theta$ is a nondegenerate bicharacter for~$G$, then 
\[R \deq \frac{\K}{|G|} \sum_{g,g' \in G} \theta(g,g') \, g \ot g' \]
is an R-matrix for the group ring $H \deq K[G]$. 
\end{prop}
\end{samepage}
\begin{pf}
The first of the axioms listed in Paragraph~\ref{QuasitriHopf} is obvious, because~$H$ is commutative and cocommutative. To verify the second axiom, we continue to write $R = \sum_{i=1}^m a_i \ot b_i$ as in Paragraph~\ref{QuasitriHopf}. We then have
\begin{align*}
\sum_{i,j=1}^m a_i \ot a_j \ot b_i b_j &= 
\frac{\K}{|G|^2} \sum_{g,g',h,h' \in G} \theta(g,g') \theta(h,h') \, g \ot h \ot g'h' \\
&= \frac{\K}{|G|^2} \sum_{g,g',h,h' \in G} \theta(g,g'h'^{-1}) \theta(h,h') \, g \ot h \ot g'
\end{align*}

Using the orthogonality relations, we can carry out the summation over~$h'$: 
\begin{align*}
\sum_{h' \in G} \theta(g,g'h'^{-1}) \theta(h,h') &= 
\sum_{h' \in G} \theta(g,g') \theta(g,h'^{-1}) \theta(h,h') \\
&= \theta(g,g') \sum_{h' \in G} \theta(hg^{-1},h')
= \delta_{g,h} |G| \theta(g,g') 
\end{align*}

In view of this equation, our previous expression becomes
\begin{align*}
\sum_{i,j=1}^m a_i \ot a_j \ot b_i b_j 
&= \frac{\K}{|G|^2} \sum_{g,g',h \in G}\delta_{g,h} |G| \theta(g,g')  \, g \ot h \ot g' \\
&= \frac{\K}{|G|} \sum_{g,g' \in G} \theta(g,g') \, g \ot g \ot g'
= \sum_{i=1}^m \Dh(a_i) \ot b_i
\end{align*}
which proves the second axiom. The verification of the third axiom is similar. Finally, we claim that 
the element~$R^{-1} = \frac{\K}{|G|} \sum_{g,g' \in G} \theta(g^{-1},g') \, g \ot g'$ is an inverse for~$R$. To see that, we form
\begin{align*}
R R^{-1} &= \frac{\K}{|G|^2} \sum_{g,g',h,h' \in G} \theta(g,g') \theta(h^{-1},h') \, gh \ot g'h' \\
&= \frac{\K}{|G|^2} \sum_{g,g',h,h' \in G} \theta(gh^{-1},g') \theta(h^{-1},h') \, g \ot g'h'
\end{align*}
Here, we can again carry out the summation over~$h$:
\begin{align*}
\sum_{h \in G} \theta(gh^{-1},g') \theta(h^{-1},h') &=  
\sum_{h \in G} \theta(g,g') \theta(h^{-1},g') \theta(h^{-1},h') \\
&= \theta(g,g') \sum_{h \in G} \theta(h^{-1}, g'h') 
=  \delta_{g',h'^{-1}} |G| \theta(g,g') 
\end{align*}

\enlargethispage{7pt}
\begin{samepage}
By resubstituting this relation back into the previous expression, we obtain
\begin{align*}
R R^{-1} &= \frac{\K}{|G|} \sum_{g,g',h' \in G} \delta_{g',h'^{-1}} \theta(g,g')  \, g \ot g'h' \\
&= \frac{\K}{|G|} \sum_{g,g' \in G}  \theta(g,g')  \, g \ot \G
= \sum_{g \in G}  \delta_{g,1} \, g \ot \G = \G \ot \G
\end{align*}
as asserted. 
\end{samepage}
\qed
\end{pf}

\section{A first example} \label{Sec:FirstExamp}
\subsection[The algebra structure of the first example]{} \label{AlgStrucFirst}
In this section, $K$ denotes a field that contains a primitive eighth root of unity~$\xi$, i.e., a nonzero element whose order in the multiplicative group is exactly~$8$. The existence of a primitive eighth root of unity forces that the characteristic of~$K$ is different from~$2$. The square~$\iota \deq \xi^2$ of~$\xi$ is a primitive fourth root of unity. We also assume that~$\zeta \in K$ is another fourth root of unity that is not required to be primitive. Clearly, $\zeta$ is a power of~$\iota$, but this fact will not play a role in the sequel.

Let $A$ be the algebra generated by two commuting variables~$x$ and~$y$ subject to the defining relations
\[x^4 = \ua \qquad \qquad y^2 = \frac{\K}{2}(\ua + \zeta x + x^2 - \zeta x^3)\]
Here and in the sequel, we suppress the index of the symbol~$\A$ and use the slightly ambiguous notation~$\ua$ instead.

The following very elementary lemma is important:
\begin{lem}
\[y^4 = 
\begin{cases} 
\ua &  \text{if} \; \; \zeta^2 = -\K \\
x^2 &  \text{if} \; \; \zeta^2 = \K
\end{cases}\]
\end{lem}
\begin{pf}
We have
\begin{align*}
y^4 &= \frac{\K}{4} \left((\ua + \zeta x)^2 +2(\ua + \zeta x)(\ua - \zeta x) x^2 + (\ua - \zeta x)^2 x^4\right)
\end{align*}
Combining the first and the last term, this becomes
\begin{align*}
y^4 &= \frac{\K}{4} \left(2 x^4 + 2 \zeta^2 x^2 + 2 (1 - \zeta^2 x^2) x^2\right) \\
&= \frac{\K}{2} \left(x^4 + \zeta^2 x^2 + x^2 - \zeta^2 x^4\right) 
= \frac{\K}{2} \left((\K - \zeta^2) x^4 +  (\K + \zeta^2) x^2\right) 
\end{align*}
which implies the assertion. 
\qed
\end{pf}

It comes at no surprise that~$A$ has dimension~$8$. It is in fact split semisimple:
\begin{prop}
As algebras, we have $A \cong K^8$.
\end{prop}
\begin{pf}
\begin{pflist}
\item
Clearly, $A$ is spanned by the elements~$\ua, x, x^2, x^3, y, xy, x^2y, x^3y$, so its dimension is at most~$8$. Since distinct characters are linearly independent by Dedekind's independence theorem (cf.~\cite{Jac3}, Sec.~I.3, Thm.~3, p.~25), our assertion will hold if we can exhibit eight distinct (one-dimensional) characters, i.e., algebra homomorphisms to the base field. For this, we need to distinguish cases that depend on the order of~$\zeta$.

\item
A character must map~$x$ to a fourth root of unity. If $\zeta$ is a primitive fourth root of unity, the following table shows the value of $\frac{\K}{2}(\ua + \zeta x + x^2 - \zeta x^3)$ depending on that choice:
\begin{center}
\renewcommand\arraystretch{1.5}
\begin{tabular}{|c|c|c|c|c|} \hline
$x$ & $\K$ & $-\K$ & $\zeta$ & $-\zeta$ \\ \hline
$\frac{\K}{2}(\ua + \zeta x + x^2 - \zeta x^3)$ & $\K$ & $\K$ & $\zeta^2$ & $\K$ \\
\hline
\end{tabular}
\end{center}

We can therefore define eight characters $\varepsilon_1$, $\varepsilon_2$, $\varepsilon_3$, $\varepsilon_4$, $\rho_1$, $\rho_2$, $\rho_3$, and~$\rho_4$ by requiring that they take the following values on the generators~$x$ and~$y$:
\begin{center}
\renewcommand\arraystretch{1.5}
\begin{tabular}{|c|c|c|c|c|c|c|c|c|} \hline
& $\varepsilon_1$ & $\varepsilon_2$ & $\varepsilon_3$ & $\varepsilon_4$ & 
$\rho_1$ & $\rho_2$ & $\rho_3$ & $\rho_4$ \\ \hline
$x$ & $\K$ & $-\K$ & $-\K$ & $\K$ & 
$\zeta$ & $-\zeta$ & $\zeta$ & $-\zeta$ \\
\hline
$y$ & $\K$ & $\K$ & $-\K$ & $-\K$ & 
$\zeta$ & $\K$ & $-\zeta$ & $-\K$ \\
\hline
\end{tabular}
\end{center}

\item
If $\zeta = \pm 1$, it is necessary to use the primitive fourth root of unity~$\iota = \xi^2$ that exists by our assumptions. If $\zeta = -1$, we can write down a similar table:
\begin{center}
\renewcommand\arraystretch{1.5}
\begin{tabular}{|c|c|c|c|c|} \hline
$x$ & $\K$ & $-\K$ & $\iota$ & $-\iota$ \\ \hline
$\frac{\K}{2}(\ua + \zeta x + x^2 - \zeta x^3)$ & $\K$ & $\K$ & $-\iota$ & $\iota$ \\
\hline
\end{tabular}
\end{center}

\begin{samepage}
Again, we can define eight characters $\varepsilon_1$, $\varepsilon_2$, $\varepsilon_3$, $\varepsilon_4$, $\rho_1$, $\rho_2$, $\rho_3$, and~$\rho_4$ by requiring that they take the following values on the generators~$x$ and~$y$:
\begin{center}
\renewcommand\arraystretch{1.5}
\begin{tabular}{|c|c|c|c|c|c|c|c|c|} \hline
& $\varepsilon_1$ & $\varepsilon_2$ & $\varepsilon_3$ & $\varepsilon_4$ & 
$\rho_1$ & $\rho_2$ & $\rho_3$ & $\rho_4$ \\ \hline
$x$ & $\K$ & $-\K$ & $-\K$ & $\K$ & 
$-\iota$ & $\iota$ & $-\iota$ & $\iota$ \\
\hline
$y$ & $\K$ & $\K$ & $-\K$ & $-\K$ & 
$\xi$ & $\iota \xi$ & $- \xi$ & $-\iota \xi$ \\
\hline
\end{tabular}
\end{center}
\end{samepage}

\item
The corresponding table in the case $\zeta = 1$ reads
\begin{center}
\renewcommand\arraystretch{1.5}
\begin{tabular}{|c|c|c|c|c|} \hline
$x$ & $\K$ & $-\K$ & $\iota$ & $-\iota$ \\ \hline
$\frac{\K}{2}(\ua + \zeta x + x^2 - \zeta x^3)$ & $\K$ & $\K$ & $\iota$ & $-\iota$ \\
\hline
\end{tabular}
\end{center}

Also in this case, we can define eight characters $\varepsilon_1$, $\varepsilon_2$, $\varepsilon_3$, $\varepsilon_4$, $\rho_1$, $\rho_2$, $\rho_3$, and~$\rho_4$ by requiring that they take the following values on the generators~$x$ and~$y$:
\begin{center}
\renewcommand\arraystretch{1.5}
\begin{tabular}{|c|c|c|c|c|c|c|c|c|} \hline
& $\varepsilon_1$ & $\varepsilon_2$ & $\varepsilon_3$ & $\varepsilon_4$ & 
$\rho_1$ & $\rho_2$ & $\rho_3$ & $\rho_4$ \\ \hline
$x$ & $\K$ & $-\K$ & $-\K$ & $\K$ & 
$-\iota$ & $\iota$ & $-\iota$ & $\iota$ \\
\hline
$y$ & $\K$ & $\K$ & $-\K$ & $-\K$ & 
$-\iota \xi$ & $\xi$ & $\iota \xi$ & $- \xi$ \\
\hline
\end{tabular}
\end{center}
In each case, we have found eight distinct characters to the base field, which we can use as the components of an algebra homomorphism from~$A$ to~$K^8$.
\qed
\end{pflist}
\end{pf}

\subsection[The module and the comodule structure of the first example]{} \label{ModComStrucFirst}
Two automorphisms of the algebra~$A$ will play a role in the sequel, which we denote by~$\phi$ and~$\phi'$. To define~$\phi$, we specify the images of the generators as 
\[\phi(x) = x^3 \qquad \qquad \phi(y) = x^3 y\]
This will yield a well-defined algebra homomorphism if we can show that~$\phi(x)$ and~$\phi(y)$ satisfy the same relations as~$x$ and~$y$. While the
relation $\phi(x)^4 = 1$ is obvious, the second relation holds by the following computation:
\[\phi(y)^2 = x^2 y^2 = \frac{\K}{2}(x^2 + \zeta x^3 + 1 - \zeta x) 
= \frac{\K}{2}(1 + \zeta \phi(x) + \phi(x)^2 - \zeta \phi(x)^3)\]
In the same way, we introduce~$\phi'$ by prescribing the images of the generators:
\[\phi'(x) = x \qquad \qquad \phi'(y) = x^2 y\]
In this case, it is obvious that~$\phi'(x)$ and~$\phi'(y)$ satisfy the same relations as~$x$ and~$y$.
\begin{lem}
The homomorphisms $\phi$ and~$\phi'$ both have order~$2$. In particular, they are bijective. Moreover, they commute.
\end{lem}
\begin{pf}
We clearly have $\phi^2(x) = x^9 = x$ and $\phi^2(y) = \phi(x)^3 \phi(y) = x^4 y = y$, so $\phi^2 = \id_A$. The corresponding verifications for~$\phi'$ are even easier. As we have 
$\phi(\phi'(x)) = x^3 = \phi'(\phi(x))$ and
\[\phi(\phi'(y)) = \phi(x)^2 \phi(y)  = x^2 (x^3 y) = \phi'(\phi(y))\]
$\phi$ and~$\phi'$ commute.
\qed
\end{pf}

It is the goal of this section to turn~$A$ into a Yetter-Drinfel'd Hopf algebra over the group ring 
$H \deq K[G]$ of the abelian group $G \deq \Z_2 \times \Z_2$. We use the notation
\[g_1 \deq (0,0) \qquad g_2 \deq (1,0) \qquad g_3 \deq (0,1) \qquad  g_4 \deq (1,1)\]
for the elements of~$G$. We introduce the action of~$G$ by requiring that $g_2$ acts by~$\phi$, $g_3$ acts by~$\phi'$, and~$g_4$ acts by
$\phi \circ \phi'$. Since~$G$ acts by algebra automorphisms, this action turns~$A$ into an $H$-module algebra, i.e., an algebra in the category of $H$-modules.

To introduce the coaction, we use the method from Paragraph~\ref{QuasitriHopf} and Paragraph~\ref{Bichar}. Our group~$G$ is elementary abelian and can therefore be considered as a vector space over the field~$\Z_2$ with two elements. Then~$g_2$ and~$g_3$ are the canonical basis vectors, so that we can define a nondegenerate symmetric bilinear form
\[\beta \coln G \times G \to \Z_2\]
by requiring that its values on the basis elements are 
\[\begin{pmatrix}
\beta(g_2,g_2) & \beta(g_2,g_3) \\ 
\beta(g_3,g_2) & \beta(g_3,g_3)
\end{pmatrix} = 
\begin{pmatrix}
* & 1 \\ 
1 & 0
\end{pmatrix} \]
where $*$ is either~$0$ or~$1$; the form will be nondegenerate in either case. From this bilinear form, we obtain a nondegenerate symmetric bicharacter~$\theta$ by defining 
\[\theta(g,g') \deq (-1)^{\beta(g,g')}\]
where for~$*$ we choose~$0$ if~$\zeta^2=1$ and~$1$ if~$\zeta^2=-1$, so that
\[\begin{pmatrix}
\theta(g_2,g_2) & \theta(g_2,g_3) \\ 
\theta(g_3,g_2) & \theta(g_3,g_3)
\end{pmatrix} = 
\begin{pmatrix}
\zeta^2 & -1 \\ 
-1 & 1
\end{pmatrix} \]

According to the discussion in Paragraph~\ref{QuasitriHopf} and Paragraph~\ref{Bichar}, $A$~now becomes a Yetter-Drinfel'd algebra if we introduce the coaction 
\[\delta \coln A \to H \ot A,~a \mapsto \frac{\K}{4} \sum_{g,g' \in G} \theta(g,g') \, g \ot g'.a\]
where we have used that our bicharacter is symmetric. We will need the explicit form of the coaction on the powers of the first generator:
\begin{prop}
The elements~$\ua$ and~$x^2$ are coinvariant. For~$x$ and~$x^3$, we have
\begin{align*}
&\delta(x) = \frac{\K}{2} (g_1 + g_3) \ot x + \frac{\K}{2} (g_1 - g_3) \ot x^3 \\
&\delta(x^3) = \frac{\K}{2} (g_1 - g_3) \ot x + \frac{\K}{2} (g_1 + g_3) \ot x^3
\end{align*}
\end{prop}
\begin{pf}
We have that $\phi(x^2) = \phi'(x^2) = x^2$ and therefore $g.x^2=x^2$ for all $g \in G$. In other words, $x^2$ as well as~$x^0=\ua$ are invariant elements. As pointed out in Paragraph~\ref{QuasitriHopf}, this implies that~$\ua$ and~$x^2$ are coinvariant. For~$x$, we have
\begin{align*}
\delta(x) &= \frac{\K}{4} \Big(\sum_{g \in G} \theta(g,g_1) \, g \ot g_1.x \Big)
+ \frac{\K}{4} \Big(\sum_{g \in G} \theta(g,g_2) \, g \ot g_2.x \Big) \\
&+ \frac{\K}{4} \Big(\sum_{g \in G} \theta(g,g_3) \, g \ot g_3.x \Big) 
+ \frac{\K}{4} \Big(\sum_{g \in G} \theta(g,g_4) \, g \ot g_4.x \Big) \\
&= \frac{\K}{4} \Big(\sum_{g \in G} (\theta(g,g_1) + \theta(g,g_3)) \, g \ot x \Big)
+ \frac{\K}{4} \Big(\sum_{g \in G} (\theta(g,g_2) + \theta(g,g_4)) \, g \ot x^3 \Big)
\end{align*}
Because we have $\theta(g,g_1) + \theta(g,g_3) = 0$ if $g=g_2$ or $g=g_4$, and we also have 
$\theta(g,g_2) + \theta(g,g_4) = 0$ if $g=g_2$ or $g=g_4$, the equation for~$x$ holds. The equation for~$x^3$ then also holds, because we have
\[\delta(x^3) = \delta(\phi(x)) = (\id_H \ot \phi)(\delta(x))\]
by the dimodule condition mentioned in Paragraph~\ref{YetMod}.
\qed
\end{pf}

\subsection[The coalgebra structure of the first example]{} \label{CoalgStrucFirst}
We now introduce a coalgebra structure on~$A$. The structures that we have already introduced on~$A$ are sufficient to form the algebra~$A \hat{\ot} A$ discussed in Paragraph~\ref{YetDrinfHopf}, and as we mentioned there, this algebra is not the canonical tensor product. Nonetheless, the elements $\ua \ot x$ and $x \ot \ua$ also commute in this algebra:
\begin{lem}
We have
\[(\ua \ot x) (x \ot \ua) = x \ot x = (x \ot \ua) (\ua \ot x)\]
\end{lem}
\begin{pf}
The second equation is a special case of the more general property
$(a \ot \ua) (\ua \ot a') = a \ot a'$. For the first equation, we use Proposition~\ref{ModComStrucFirst} and the fact that $g_3.x = \phi'(x) = x$ to get
\begin{align*}
(\ua \ot x) (x \ot \ua) = \frac{\K}{2} (g_1 + g_3).x \ot x + \frac{\K}{2} (g_1 - g_3).x \ot x^3
= x \ot x
\end{align*}
as asserted.
\qed
\end{pf}

The construction of the coproduct~$\da$ can now be based on the following proposition:
\begin{prop}
There is a unique algebra homomorphism $\da \coln A \to A \hat{\ot} A$ with the property that
\[\da(x) =  \frac{\K}{2} (x \ot x + x \ot x^3 + x^3 \ot x - x^3 \ot x^3)
\qquad \da(y) = y \ot y\]
\end{prop}
\begin{pf}
\begin{pflist}
\item
We need to show that the potential images 
\[x' \deq  \frac{\K}{2} (x \ot x + x \ot x^3 + x^3 \ot x - x^3 \ot x^3)
\qquad \text{and} \qquad y' \deq y \ot y\]
satisfy the defining relations of the algebra~$A$. Using the preceding lemma, it is not very difficult to verify the equation $x'^2 = x^2 \ot x^2$, so that $x'^4 = \ua \ot \ua$.

\item
Next, we compute the square of~$y'$ in $A \hat{\ot} A$. First, we note that the relation
$y^2 = \frac{\K}{2}((\ua + x^2) + \zeta x(\ua - x^2))$ implies
\[(\ua + x^2) y^2 = \ua + x^2 \qquad \text{and} \qquad 
(\ua - x^2) y^2 = \zeta x(\ua - x^2)\]
We have
\begin{align*}
&y'^2 = \frac{\K}{4} \sum_{i,j=1}^4 \theta(g_i,g_j) y (g_i.y) \ot (g_j.y) y 
\end{align*}
We look at the terms with~$i$ fixed individually, and carry out the summation over~$j$. For $i=1$, we have
\begin{align*}
\sum_{j=1}^4 \theta(g_1,g_j) (g_j.y) y &= \sum_{j=1}^4 (g_j.y) y = 
y^2 + x^3 y^2 + x^2 y^2 + x y^2 = \ua + x + x^2 + x^3   
\end{align*}

For $i=2$, we have
\begin{align*}
&\sum_{j=1}^4 \theta(g_2,g_j) (g_j.y) y = 
y^2 + \zeta^2 x^3 y^2 - x^2 y^2 - \zeta^2 x y^2 \\
&= (\ua - x^2) y^2 - \zeta^2 x (\ua - x^2) y^2 
= \zeta x(\ua - x^2) - \zeta^3 x^2 (\ua - x^2) \\
&= \zeta^3 \ua + \zeta x - \zeta^3 x^2 - \zeta x^3
\end{align*}

For $i=3$, we have
\begin{align*}
&\sum_{j=1}^4 \theta(g_3,g_j) (g_j.y) y = 
y^2 - x^3 y^2 + x^2 y^2 - x y^2 = \ua - x  + x^2 - x^3
\end{align*}

For $i=4$, we have
\begin{align*}
&\sum_{j=1}^4 \theta(g_4,g_j) (g_j.y) y = 
y^2 - \zeta^2 x^3 y^2 - x^2 y^2 + \zeta^2 x y^2 \\
&= (\ua - x^2) y^2 + \zeta^2 x (\ua - x^2) y^2 
= \zeta x(\ua - x^2) + \zeta^3 x^2 (\ua - x^2) \\
&= - \zeta^3 \ua + \zeta x + \zeta^3 x^2 - \zeta x^3 
\end{align*}

In total, we obtain
\begin{align*}
4 y'^2 &= y^2 \ot (\ua + x + x^2 + x^3) +
x^3 y^2 \ot (\zeta^3 \ua + \zeta x - \zeta^3 x^2 - \zeta x^3)  \\
&+ x^2 y^2 \ot (\ua - x  + x^2 - x^3)  +
x y^2 \ot (- \zeta^3 \ua + \zeta x + \zeta^3 x^2 - \zeta x^3 )  
\end{align*}

Collecting terms, this becomes
\begin{align*}
4 y'^2 &= (\ua + \zeta^3 x^3 + x^2 - \zeta^3 x) y^2 \ot \ua +
(\ua + \zeta x^3 - x^2 + \zeta x) y^2 \ot x \\
&+ (\ua - \zeta^3 x^3 + x^2 + \zeta^3 x) y^2 \ot x^2   +
(\ua - \zeta x^3 - x^2 - \zeta x) y^2 \ot x^3 
\end{align*}
or, with the help of the equations recorded above,
\begin{align*}
4 y'^2 &= (\ua + \zeta^4 x^4 + x^2 - \zeta^4 x^2) \ot \ua +
(\zeta x + \zeta x^3 - \zeta x^3 + \zeta x) \ot x \\
&+ (\ua - \zeta^4 x^4 + x^2 + \zeta^4 x^2) \ot x^2   +
(\zeta x - \zeta x^3 - \zeta x^3 - \zeta x) \ot x^3 
\end{align*}
which simplifies to
\begin{align*}
y'^2 &= \frac{\K}{2} (\ua \ot \ua + \zeta x \ot x + x^2  \ot x^2 - \zeta x^3 \ot x^3) 
\end{align*}

\enlargethispage{8pt}
\begin{samepage}
On the other hand, we have
\begin{align*}
\frac{\K}{2}(&1 +\zeta x' + x'^2 - \zeta x'^3) = \\
&\frac{\K}{2}(\ua \ot \ua + \frac{\zeta}{2} (x \ot x + x \ot x^3 + x^3 \ot x - x^3 \ot x^3) \\
&\mspace{20mu}
+ x^2 \ot x^2 - \frac{\zeta}{2} (x^3 \ot x^3 + x^3 \ot x + x \ot x^3 - x \ot x)) \\
&= 
\frac{\K}{2}(\ua \ot \ua + x^2 \ot x^2 + \zeta (x \ot x - x^3 \ot x^3))
\end{align*}
which establishes that $y'^2 = \frac{\K}{2}(1 +\zeta x' + x'^2 - \zeta x'^3)$.
\end{samepage}

\item
We claim that $(y \ot y)(x \ot x) = (x \ot x)(y \ot y)$. We have by Proposition~\ref{ModComStrucFirst} that
\begin{align*}
(\ua \ot x)(y \ot \ua) &= (x^{\1}.y) \ot x^{\2}
= \frac{\K}{2} (g_1 + g_3).y \ot x + \frac{\K}{2} (g_1 - g_3).y \ot x^3 \\
&= \frac{\K}{2} (\ua + x^2) y \ot x + \frac{\K}{2} (\ua - x^2) y \ot x^3
\end{align*}
This implies that
\begin{align*}
(x \ot x)(y \ot y) 
&= \frac{\K}{2} (x + x^3) y \ot x y + \frac{\K}{2} (x - x^3) y \ot x^3 y \\
&= \frac{\K}{2} (xy  \ot xy + x^3y \ot xy + xy \ot x^3 y - x^3 y \ot x^3 y)
\end{align*}

On the other hand, we get from the symmetry of~$\theta$ that
\begin{align*}
(\ua \ot y)(x \ot \ua) &= (y^{\1}.x) \ot y^{\2} \\
&= \frac{\K}{4} \sum_{i,j=1}^4 \theta(g_i,g_j)  (g_i.x) \ot (g_j.y) 
= x^{\2} \ot x^{\1}.y
\end{align*}
which by the calculation just above yields
\begin{align*}
(\ua \ot y)(x \ot \ua) &= \frac{\K}{2} x \ot (\ua + x^2) y + \frac{\K}{2} x^3 \ot (\ua - x^2)y
\end{align*}
This implies that
\begin{align*}
(y \ot y)(x \ot x) &= \frac{\K}{2} xy \ot (x + x^3) y + \frac{\K}{2} x^3 y \ot (x - x^3)y \\
&= \frac{\K}{2} (xy \ot xy + xy \ot x^3 y +  x^3 y \ot xy - x^3 y \ot x^3 y)
\end{align*}
proving our claim.

\item
We have seen in Proposition~\ref{ModComStrucFirst} that~$x^2$ is invariant and coinvariant. This implies that the elements $x^2 \ot \ua$ and $\ua \ot x^2$ are central in~$A \hat{\ot} A$.
From the lemma proved at the beginning of this paragraph, we get
$(x \ot x)^2 = x^2 \ot x^2$, which then implies by the centrality just mentioned that 
$(x \ot x)^3 = x^3 \ot x^3$. The claim just established in the preceding step therefore implies that 
\[(y \ot y)(x^3 \ot x^3) = (x^3 \ot x^3)(y \ot y)\]
and in view of the centrality of~$x^2 \ot \ua$ and $\ua \ot x^2$ also
\[(y \ot y)(x^3 \ot x) = (x^3 \ot x)(y \ot y) \quad \text{and} \quad 
(y \ot y)(x \ot x^3) = (x \ot x^3)(y \ot y)\]
All these relations combined yield that~$x'y' = y'x'$, which was the last relation for the elements~$x'$ and~$y'$ that we needed to establish. Our proposition is therefore proved.
\qed
\end{pflist}
\end{pf}

We now define the coproduct to be the algebra homomorphism $\da \coln A \to A \hat{\ot} A$ uniquely determined by the properties stated in the preceding proposition. We claim that
$(\phi \ot \phi) \circ \da = \da \circ \phi$. Since both the left and the right-hand side of this equation are algebra homomorphisms, it suffices to check this on the generators~$x$ and~$y$. For~$x$, we have on the one hand
\begin{align*}
(\phi \ot \phi) (\da(x)) &= 
\frac{\K}{2}(\phi \ot \phi) (x \ot x + x \ot x^3 + x^3 \ot x - x^3 \ot x^3) \\
&= \frac{\K}{2} (x^3 \ot x^3 + x^3 \ot x + x \ot x^3 - x \ot x) 
\end{align*}

On the other hand, we have computed the square of~$x'$ in the preceding proof and therefore get
\begin{align*}
\da(\phi(x)) &= \da(x^3) = \da(x)^3 = \da(x)^2 \da(x) = (x^2 \ot x^2) \da(x)\\
&=
\frac{\K}{2}(x^2 \ot x^2) (x \ot x + x \ot x^3 + x^3 \ot x - x^3 \ot x^3) \\
&= \frac{\K}{2} (x^3 \ot x^3 + x^3 \ot x + x \ot x^3 - x \ot x) 
\end{align*}

For~$y$, we can use the formula for $(x \ot x) (y \ot y)$ found in the preceding proof to obtain
\begin{align*}
&\da(\phi(y)) = \da(x^3 y) = \da(x)^3 \da(y) \\
&= \frac{\K}{2} (x^3 \ot x^3 + x^3 \ot x + x \ot x^3 - x \ot x) (y \ot y) \\
&= \frac{\K}{2} (x^2 \ot x^2 + x^2 \ot \ua + \ua \ot x^2 - \ua \ot \ua) (x \ot x) (y \ot y) \\
&= \frac{\K}{4} (x^2 \ot x^2 + x^2 \ot \ua + \ua \ot x^2 - \ua \ot \ua) \\
&\mspace{40mu}
(xy  \ot xy + x^3y \ot xy + xy \ot x^3 y - x^3 y \ot x^3 y)\
\end{align*}

Expanding the terms, we find
\begin{align*}
\da(\phi(y)) &= 
\frac{\K}{4} (x^3y \ot x^3y + xy \ot x^3y + x^3y \ot x y - xy \ot xy)  \\
&+ \frac{\K}{4} (x^3 y \ot xy + xy \ot xy + x^3y \ot x^3 y - xy \ot x^3 y) \\
&+ \frac{\K}{4} (xy \ot x^3 y + x^3y \ot x^3 y + xy \ot x y - x^3 y \ot x y) \\
&- \frac{\K}{4} (xy \ot xy + x^3y \ot xy + xy \ot x^3 y - x^3 y \ot x^3 y)
\end{align*}
which reduces to $\da(\phi(y)) = x^3y \ot x^3y$. On the other hand, we have
\begin{align*}
(\phi \ot \phi)(\da(y)) &= (\phi \ot \phi)(y \ot y) = x^3 y \ot x^3 y
\end{align*}
A similar, but simpler computation shows that $(\phi' \ot \phi') \circ \da = \da \circ \phi'$. Because $A \hat{\ot} A$ carries the diagonal module structure, this yields
$\da(g.a) = g.\da(a)$ for all~$g \in G$ and all~$a \in A$. In other words, $\da$ is $H$-linear. As mentioned in Paragraph~\ref{QuasitriHopf}, it is then automatically $H$-colinear.

We have seen in Paragraph~\ref{AlgStrucFirst} that, regardless of the value of~$\zeta$, there is an algebra homomorphism~$\varepsilon_1 \coln A \to K$ which takes the values
\[\varepsilon_1(x) = \K = \varepsilon_1(y)\]
on the generators. We define the counit to be this algebra homomorphism; i.e., we set~$\ea \deq \varepsilon_1$. It is immediate to check that
$\ea \circ \phi = \ea = \ea \circ \phi'$, which implies that \mbox{$\ea(g.a) = \ea(a)$} for 
all~$g \in G$ and all~$a \in A$, so that $\ea$ is $H$-linear. As just discussed in the case of the coproduct, $\ea$~is then automatically also $H$-colinear.

We claim that $(\ea \ot \id_A)(\da(a)) = \K \ot a$ for all~$a \in A$. Since both the left and the right-hand side depend multiplicatively on~$a$, it suffices to check this for the generators~$x$ and~$y$. For~$x$, we have
\begin{align*}
(\ea \ot \id_A) (\da(x)) &= 
\frac{\K}{2}(\ea \ot \id_A) (x \ot x + x \ot x^3 + x^3 \ot x - x^3 \ot x^3) \\
&= \frac{\K}{2} (\K \ot x + \K \ot x^3 + \K \ot x - \K \ot x^3) = \K \ot x
\end{align*}
For~$y$, we have $(\ea \ot \id_A) (\da(y)) = (\ea \ot \id_A) (y \ot y) = \K \ot y$. The verification of the equation $(\id_A \ot \, \ea)(\da(a)) = a \ot \K$ is very similar.

\subsection[A basis consisting of group-like elements for the first example]{} \label{BasisGrouplFirst}
It remains to be shown that the coproduct is coassociative. This can be done by verifying that the two algebra homomorphisms~$(\da \ot \id_A) \circ \da$ and $(\id_A \ot \, \da) \circ \da$ agree on the generators~$x$ and~$y$, up to the canonical isomorphism between the algebras
$(A \hat{\ot} A) \hat{\ot} A$ and $A \hat{\ot} (A \hat{\ot} A)$. However, we proceed differently and exhibit a basis consisting of group-like elements, which also proves coassociativity:
\begin{defn}
We define $\omega_1 \deq \ua$,
\[\omega_2 \deq \frac{\K}{2}(1 + \iota \zeta^2) x + \frac{\K}{2}(1 - \iota \zeta^2) x^3 \qquad
\omega_3 \deq \frac{\K}{2}(1 - \iota \zeta^2) x + \frac{\K}{2}(1 + \iota \zeta^2) x^3\] 
and $\omega_4 \deq x^2$. Furthermore, we define
\[\eta_1 \deq y \qquad \eta_2 \deq x^3 y \qquad \eta_3 \deq x^2 y \qquad \eta_4 \deq x y\]
\end{defn}

The following proposition shows that this definition is appropriate:
\begin{prop}
The elements $\omega_1, \omega_2, \omega_3, \omega_4, \eta_1, \eta_2, \eta_3, \eta_4$ form a basis of~$A$. They are group-like. 
\end{prop}
\begin{pf}
\begin{pflist}
\item
By Proposition~\ref{AlgStrucFirst}, the elements~$\ua, x, x^2, x^3, y, xy, x^2y, x^3y$ form a basis of~$A$. Since the matrix
\[\begin{pmatrix}
1 + \iota \zeta^2 & 1 - \iota \zeta^2 \\
1 - \iota \zeta^2 & 1 + \iota \zeta^2
\end{pmatrix}\] 
has determinant $4 \iota \zeta^2 \neq 0$, the new elements also form a basis. 

\item
It is obvious that $\omega_1 \deq \ua$ is group-like, and we have already seen in the first step of the proof of Proposition~\ref{CoalgStrucFirst} that $\omega_4 \deq x^2$ is group-like. We have
\[\da(x) =  \frac{\K}{2} (x \ot x + x \ot x^3 + x^3 \ot x - x^3 \ot x^3)\]
by definition, and in the discussion of the $H$-linearity of~$\da$, we have seen that
\begin{align*}
\da(x^3) &= \frac{\K}{2} (x^3 \ot x^3 + x^3 \ot x + x \ot x^3 - x \ot x) 
\end{align*}

Therefore we have
\begin{align*}
\da(&\omega_2) = 
\frac{\K}{2}(1 + \iota \zeta^2) \da(x) + \frac{\K}{2}(1 - \iota \zeta^2) \da(x^3) \\
&= \frac{\K}{2} (\iota \zeta^2 x \ot x + x \ot x^3 + x^3 \ot x - \iota \zeta^2 x^3 \ot x^3) \\
&= \frac{\K}{4}\Big((1 + \iota \zeta^2) x + (1 - \iota \zeta^2) x^3\Big) \ot
\Big((1 + \iota \zeta^2) x + (1 - \iota \zeta^2) x^3\Big) = \omega_2 \ot \omega_2
\end{align*}

As we will discuss in detail in Paragraph~\ref{ProdGroupl}, the product of two group-like elements is in general not again group-like. But the situation is better for group-like elements that are invariant or coinvariant (cf.~\cite{SoDiss}, Prop.~1.5.1, p.~13; \cite{SoYp}, Prop.~1.5.1, p.~12), and since~$\omega_4$ has these two properties, we get 
that~$\omega_3 = \omega_2 \omega_4$ is also group-like. 

\item
The element~$\eta_1 = y$ is group-like by construction. We have already seen that~$\eta_2 = x^3 y$ is group-like when we discussed the $H$-linearity of the coproduct at the end of Paragraph~\ref{CoalgStrucFirst}. But then~$\eta_3 = \omega_4 \eta_1$ and~$\eta_4 = \omega_4 \eta_2$ are group-like for the reasons just cited. 
\qed
\end{pflist}
\end{pf}

To exhibit the essential features of the algebra~$A$ in Section~\ref{Sec:Cores}, it is important to have the products of the group-like elements explicitly available. It is not complicated to check that~$\omega_2^2 = \omega_1 = \ua$, and it is immediate that $\omega_3 = \omega_2 \omega_4$. Based on this information, it is not too complicated to verify the following table for the products~$\omega_i \omega_j$:
\begin{center}
\renewcommand\arraystretch{1.5}
\begin{tabular}{|c|c|c|c|c|}
\hline 
& $\omega_1$ & $\omega_2$ & $\omega_3$ & $\omega_4$   \\ 
\hline 
$\omega_1$ & $\omega_1$ & $\omega_2$ & $\omega_3$ & $\omega_4$   \\ 
\hline 
$\omega_2$ & $\omega_2$ & $\omega_1$ & $\omega_4$ & $\omega_3$ \\ 
\hline 
$\omega_3$ & $\omega_3$ & $\omega_4$ & $\omega_1$ & $\omega_2$ \\ 
\hline 
$\omega_4$ & $\omega_4$ & $\omega_3$ & $\omega_2$ & $\omega_1$ \\ 
\hline 
\end{tabular} 
\end{center}
Note that we do not need to say which factor comes first, because~$A$ is commutative. This table shows that
$\Span(\omega_1, \omega_2, \omega_3, \omega_4)$ is a subalgebra isomorphic to~$K[G]$ under the isomorphism $g_i \mapsto \omega_i$.

The table for the products $\omega_i \eta_j = \eta_j \omega_i$ is easily verified from the definitions:
\begin{center}
\renewcommand\arraystretch{1.5}
\begin{tabular}{|c|c|c|c|c|}
\hline 
& $\omega_1$ & $\omega_2$ & $\omega_3$ & $\omega_4$   \\ 
\hline 
$\eta_1$ & $\eta_1$ & 
$\frac{\K}{2}(1 - \iota \zeta^2) \eta_2 + \frac{\K}{2} (1 + \iota \zeta^2) \eta_4$ & $\frac{\K}{2}(1 + \iota \zeta^2) \eta_2 + \frac{\K}{2} (1 - \iota \zeta^2) \eta_4$ & 
$\eta_3$ \\ 
\hline 
$\eta_2$ & $\eta_2$ & 
$\frac{\K}{2}(1 + \iota \zeta^2) \eta_1 + \frac{\K}{2} (1 - \iota \zeta^2) \eta_3$ & 
$\frac{\K}{2}(1 - \iota \zeta^2) \eta_1 + \frac{\K}{2} (1 + \iota \zeta^2) \eta_3$ & $\eta_4$ \\ 
\hline 
$\eta_3$ & $\eta_3$ & 
$\frac{\K}{2}(1 + \iota \zeta^2) \eta_2 + \frac{\K}{2} (1 - \iota \zeta^2) \eta_4$ & $\frac{\K}{2}(1 - \iota \zeta^2) \eta_2 + \frac{\K}{2} (1 + \iota \zeta^2) \eta_4$ & $\eta_1$ \\ 
\hline 
$\eta_4$ & $\eta_4$ & 
$\frac{\K}{2}(1 - \iota \zeta^2) \eta_1 + \frac{\K}{2} (1 + \iota \zeta^2) \eta_3$ & 
$\frac{\K}{2}(1 + \iota \zeta^2) \eta_1 + \frac{\K}{2} (1 - \iota \zeta^2) \eta_3$ & $\eta_2$ \\ 
\hline 
\end{tabular} 
\end{center}

For the products $\eta_i \eta_j = \eta_j \eta_i$, we break the table into two parts. The first part lists the products that involve~$\eta_1$ and~$\eta_2$:
\begin{center}
\renewcommand\arraystretch{1.5}
\begin{tabular}{|c|c|c|}
\hline  
 & $\eta_1$ & $\eta_2$ \\ 
\hline 
$\eta_1$ & 
$\frac{\K}{2} \omega_1 - \frac{\iota}{2\zeta} \omega_2 
+ \frac{\iota}{2\zeta} \omega_3 + \frac{\K}{2} \omega_4$ & 
$\frac{\zeta}{2} \omega_1 + \frac{\K}{2} \omega_2 
+ \frac{\K}{2} \omega_3 - \frac{\zeta}{2} \omega_4$ \\ 
\hline 
$\eta_2$ & 
$\frac{\zeta}{2} \omega_1 + \frac{\K}{2} \omega_2 
+ \frac{\K}{2} \omega_3 - \frac{\zeta}{2} \omega_4$ & 
$\frac{\K}{2} \omega_1 + \frac{\iota}{2\zeta} \omega_2 
- \frac{\iota}{2\zeta} \omega_3 + \frac{\K}{2} \omega_4$ \\ 
\hline 
$\eta_3$ & 
$\frac{\K}{2} \omega_1 + \frac{\iota}{2\zeta} \omega_2 
- \frac{\iota}{2\zeta} \omega_3 + \frac{\K}{2} \omega_4$ & 
$- \frac{\zeta}{2} \omega_1 + \frac{\K}{2} \omega_2 
+ \frac{\K}{2} \omega_3 + \frac{\zeta}{2} \omega_4$  \\ 
\hline 
$\eta_4$ & 
$- \frac{\zeta}{2} \omega_1 + \frac{\K}{2} \omega_2 
+ \frac{\K}{2} \omega_3 + \frac{\zeta}{2} \omega_4$ & 
$\frac{\K}{2} \omega_1 - \frac{\iota}{2\zeta} \omega_2 
+ \frac{\iota}{2\zeta} \omega_3 + \frac{\K}{2} \omega_4$  \\ 
\hline 
\end{tabular} 
\end{center}

The first entry in this table arises from the facts that $\eta_1 = y$ and
\[- \frac{\iota}{\zeta} \omega_2 + \frac{\iota}{\zeta} \omega_3 = \zeta x - \zeta x^3 \]
The remaining entries follow comparatively easily from the first one.

The second part lists the products that do not involve~$\eta_1$ and~$\eta_2$:
\begin{center}
\renewcommand\arraystretch{1.5}
\begin{tabular}{|c|c|c|c|c|}
\hline  
 & $\eta_3$ & $\eta_4$ \\ 
\hline 
$\eta_3$ & 
$\frac{\K}{2} \omega_1 - \frac{\iota}{2\zeta} \omega_2 
+ \frac{\iota}{2\zeta} \omega_3 + \frac{\K}{2} \omega_4$ & 
$\frac{\zeta}{2} \omega_1 + \frac{\K}{2} \omega_2 
+ \frac{\K}{2} \omega_3 - \frac{\zeta}{2} \omega_4$ \\ 
\hline 
$\eta_4$  & 
$\frac{\zeta}{2} \omega_1 + \frac{\K}{2} \omega_2 
+ \frac{\K}{2} \omega_3 - \frac{\zeta}{2} \omega_4$  & 
$\frac{\K}{2} \omega_1 + \frac{\iota}{2\zeta} \omega_2 
- \frac{\iota}{2\zeta} \omega_3 + \frac{\K}{2} \omega_4$ \\ 
\hline 
\end{tabular} 
\end{center}
The entries in the second part follow again relatively easily from the entries in the first part.

Clearly, the action of the group~$G$ on~$A$ permutes the group-like elements. The following lemma describes how:
\begin{lem}
\leavevmode
\begin{thmlist}
\item 
$\omega_1$ and~$\omega_4$ are fixed points of the $G$-action. Moreover, we have $g_3.\omega_i = \omega_i$ for all $i=1,2,3,4$, while $g_2.\omega_2 = \omega_3$.

\item
$G$ acts simply transitively on the set~$\{\eta_1, \eta_2, \eta_3, \eta_4\}$: If we define
$\eta_{g_i} \deq \eta_i$, we have $g_i.\eta_{g_j} = \eta_{g_i g_j}$.
\end{thmlist}
\end{lem}
\begin{pf}
The first assertion about the~$\omega_i$ follows directly from the construction of the action in Paragraph~\ref{ModComStrucFirst}. The equation $g_i.\eta_{g_j} = \eta_{g_i g_j}$ can be checked on generators, i.e., for~$i=2$ and~$i=3$, where it again follows directly from the definitions.
\qed
\end{pf}

\subsection[The antipode in the first example]{} \label{AntipFirst}
The remaining structure element that we need to construct is the antipode. We define it to be the
unique $K$-linear endomorphism~$\sa \coln A \to A$ that is given on the basis from Proposition~\ref{BasisGrouplFirst} by the formulas 
\begin{align*}
\sa(\omega_1) &= \omega_1 \qquad \sa(\omega_4) = \omega_4 \quad &
\sa(\omega_2) &= \omega_2 \qquad \sa(\omega_3) = \omega_3 \\
\sa(\eta_1) &= 
\frac{\K}{2}\left(\eta_1 + \frac{\K}{\zeta} \eta_2 + \eta_3 - \frac{\K}{\zeta} \eta_4 \right) \quad &
\sa(\eta_2) &= 
\frac{\K}{2}\left(\frac{\K}{\zeta}\eta_1 + \eta_2 - \frac{\K}{\zeta} \eta_3 + \eta_4 \right)  \\
\sa(\eta_3) &= 
\frac{\K}{2}\left(\eta_1 - \frac{\K}{\zeta} \eta_2 + \eta_3 + \frac{\K}{\zeta} \eta_4 \right) \quad &
\sa(\eta_4) &= 
\frac{\K}{2}\left(- \frac{\K}{\zeta}\eta_1 +  \eta_2 + \frac{\K}{\zeta} \eta_3 + \eta_4 \right)
\end{align*}

We will now prove that this endomorphism has the required properties:
\begin{prop}
$\sa$ is an antipode for~$A$.
\end{prop}
\begin{pf}
The defining equation  
\[\sa(a_\1) a_\2 = a_\1 \sa(a_\2) = \ea(a) \ua \]
of the antipode is $K$-linear, so that it will be satisfied if it is satisfied on the basis used for its definition. By Proposition~\ref{BasisGrouplFirst}, this basis consists of group-like elements, for which the equation just says that the image of a group-like element under the antipode is its inverse. For the elements $\omega_1$, $\omega_2$, $\omega_3$, and~$\omega_4$, this follows from their multiplication table given in Paragraph~\ref{BasisGrouplFirst}. For the basis element~$\eta_1$, we have by another multiplication table that
\begin{align*}
\eta_1 \sa(\eta_1) &= 
\frac{\K}{2}(\eta_1 \eta_1 + \frac{\K}{\zeta} \eta_1 \eta_2 
+ \eta_1 \eta_3 - \frac{\K}{\zeta} \eta_1 \eta_4) \\
&= \frac{\K}{4}((\omega_1 - \frac{\iota}{\zeta} \omega_2 
+ \frac{\iota}{\zeta} \omega_3 + \omega_4) 
+ \frac{\K}{\zeta} (\zeta \omega_1 +  \omega_2 +  \omega_3 - \zeta \omega_4) \\
&+ (\omega_1 + \frac{\iota}{\zeta} \omega_2 
- \frac{\iota}{\zeta} \omega_3 + \omega_4) 
- \frac{\K}{\zeta} (- \zeta \omega_1 + \omega_2 +  \omega_3 + \zeta \omega_4))
= \omega_1
\end{align*}

It would be possible to treat the remaining basis elements~$\eta_2$, $\eta_3$, and~$\eta_4$ in the same way, and then to use Lemma~\ref{YetDrinfHopf} to conclude that~$\sa$ is $H$-linear and colinear. However, it is easier to note that the definition implies immediately that~$\sa$ commutes with~$\phi$ and~$\phi'$, i.e., is $H$-linear, and therefore also colinear, as mentioned in Paragraph~\ref{QuasitriHopf}. From the equation $\eta_1 \sa(\eta_1)= \omega_1$, we then get by applying~$\phi$ that
$\eta_2 \sa(\eta_2)= \omega_1$, by applying~$\phi'$ that
$\eta_3 \sa(\eta_3)= \omega_1$, and by applying~$\phi \circ \phi'$ that
$\eta_4 \sa(\eta_4)= \omega_1$. Since $A$ is commutative, this is sufficient.
\qed
\end{pf}

\section{A second example} \label{Sec:SecExamp}
\subsection[The algebra of the second example]{} \label{AlgSec}
In this section, $K$ denotes a field that contains a primitive fourth root of unity~$\iota$. As in~Section~\ref{Sec:FirstExamp}, the existence of a primitive fourth root of unity forces that the characteristic of~$K$ is different from~$2$. Also as in Section~\ref{Sec:FirstExamp}, we assume that~$\zeta \in K$ is another fourth root of unity that is not required to be primitive. 

We introduce a new algebra, also denoted by~$A$, as the algebra generated by two noncommuting variables~$x$ and~$y$ subject to the defining relations
\[x^4 = \ua \qquad \qquad xy = yx^3 \qquad \qquad 
y^2 = \frac{\K}{2}(\zeta \ua + x - \zeta x^2 + x^3)\]

As in Paragraph~\ref{AlgStrucFirst}, the fourth power of~$y$ depends on the question whether or not~$\zeta$ is primitive:
\begin{lem}
\[y^4 = 
\begin{cases} 
\ua &  \text{if} \; \; \zeta^2 = \K \\
x^2 &  \text{if} \; \; \zeta^2 = -\K
\end{cases}\]
\end{lem}
This fact can be shown as in Lemma~\ref{AlgStrucFirst}, or in fact be deduced from there.

Because~$A$ is now noncommutative, we need a different argument to establish its dimension: 
\begin{prop}
The elements $\ua, x, x^2, x^3, y, xy, x^2y, x^3y$ form a basis of~$A$. In particular, $A$ has dimension~$8$. 
\end{prop}
\begin{pf}
It is clear that these elements span~$A$ linearly over~$K$. To prove linear independence, we consider two $8 \times 8$-matrices~$X$ and~$Y$ over~$K$, which we construct as $2 \times 2$-blocks of $4 \times 4$-matrices: We set
\[X \deq \begin{pmatrix} X_1 & 0 \\ 0 & X_2 \end{pmatrix}
\qquad  \text{and} \qquad 
Y \deq \begin{pmatrix} 0 & Y_1 \\ Y_2 & 0 \end{pmatrix}\]

\begin{samepage}
where the~zeroes denote $4 \times 4$-zero matrices and the other blocks are the $4 \times 4$-matrices
\[X_1 = X_2 \deq \begin{pmatrix}
0 & 0 & 0 & 1 \\
1 & 0 & 0 & 0 \\
0 & 1 & 0 & 0 \\
0 & 0 & 1 & 0  
\end{pmatrix}\]
for~$X$ and
\[Y_1 \deq \frac{\K}{2} \begin{pmatrix}
\zeta & 1 & -\zeta & 1 \\
1 & -\zeta & 1 & \zeta \\
- \zeta & 1 & \zeta & 1 \\
1 & \zeta & 1 & -\zeta
\end{pmatrix}
\qquad  \text{and} \qquad
Y_2 \deq
\begin{pmatrix}
1 & 0 & 0 & 0 \\
0 & 0 & 0 & 1 \\
0 & 0 & 1 & 0 \\
0 & 1 & 0 & 0 
\end{pmatrix}\]
for~$Y$.
\end{samepage}

These matrices satisfy the above relations: It is obvious that 
$X^4 = E_8$, the $8 \times 8$-identity matrix. The relation $XY = YX^3$ follows from the identities
\[X_1 Y_1 = \frac{\K}{2} \begin{pmatrix}
1 & \zeta & 1 & -\zeta \\
\zeta & 1 & -\zeta & 1 \\
1 & -\zeta & 1 & \zeta \\
- \zeta & 1 & \zeta & 1
\end{pmatrix} = Y_1 X_2^3 \]
and 
\[X_2 Y_2 = 
\begin{pmatrix}
0 & 1 & 0 & 0 \\
1 & 0 & 0 & 0 \\
0 & 0 & 0 & 1 \\
0 & 0 & 1 & 0 
\end{pmatrix}= Y_2 X_1^3 \]

The relation $Y^2 = \frac{\K}{2}(\zeta E_8 + X - \zeta X^2 + X^3)$ means for the $4 \times 4$-blocks that
$Y_1 Y_2 = Y_2 Y_1 = \frac{\K}{2}(\zeta E_4 + X_1 - \zeta X_1^2 + X_1^3)$. Now~$Y_1$ and~$Y_2$ indeed commute, since
\[Y_1 Y_2 = \frac{\K}{2} \begin{pmatrix}
\zeta & 1 & -\zeta & 1 \\
1 & \zeta & 1 & -\zeta \\
- \zeta & 1 & \zeta & 1 \\
1 & -\zeta & 1 & \zeta
\end{pmatrix}
= Y_2 Y_1\]
and this matrix is also equal to~$\frac{\K}{2}(\zeta E_4 + X_1 - \zeta X_1^2 + X_1^3)$. This implies that there is an algebra homomorphism from~$A$ to~$M(8 \times 8,K)$
that maps~$x$ to~$X$ and $y$~to~$Y$, so that $E_8, X, X^2, X^3, Y, XY, X^2Y, X^3Y$
would be linearly dependent if $\ua, x, x^2, x^3, y, xy, x^2y, x^3y$ would be linearly dependent. But these matrices are linearly independent, as follows from the linear independence of~$E_4, X_1, X_1^2, X_1^3$ in the left upper $4 \times 4$-block and the linear independence of~$Y_2, X_1 Y_2, X_1^2 Y_2, X_1^3 Y_2$ in the left lower $4 \times 4$-block.
\qed
\end{pf}

It should be noted that, after proving this proposition, we can interpret the matrices~$X$ and~$Y$ as the matrix representations of the left multiplications by~$x$ and~$y$ with respect to the basis now established.

\subsection[The module and the comodule structure of the second example]{} \label{ModComStrucSec}
In almost the same way as in Paragraph~\ref{ModComStrucFirst}, we turn~$A$ into a Yetter-Drinfel'd algebra over the group ring $H = K[G]$ of the group $G = \Z_2 \times \Z_2$. We begin by defining two algebra homomorphisms~$\phi$ and~$\phi'$ by the same formulas as in the previous case: First,~$\phi$ is defined on generators by 
\[\phi(x) = x^3 \qquad \qquad \phi(y) = x^3 y\]
Again, this will yield a well-defined algebra homomorphism if we can show that~$\phi(x)$ and~$\phi(y)$ satisfy the same relations as~$x$ and~$y$. While the relations 
\[\phi(x)^4 = \ua \qquad \text{and} \qquad \phi(x) \phi(y) = \phi(y) \phi(x)^3\]
are more or less obvious, the last relation holds by the following computation:
\begin{align*}
\phi(y)^2 &= (x^3 y) (x^3 y) = x^3 (x y) y = y^2 
= \frac{\K}{2}(\zeta \ua + x - \zeta x^2 + x^3) \\
&= \frac{\K}{2}(\zeta \ua + \phi(x) - \zeta \phi(x)^2 + \phi(x)^3)
\end{align*}
We also introduce the automorphism~$\phi'$ by the same formulas as in Paragraph~\ref{ModComStrucFirst} and require that
\[\phi'(x) = x \qquad \qquad \phi'(y) = x^2 y\]
Again, it is very easy to verify that~$\phi'(x)$ and~$\phi'(y)$ satisfy the same relations as~$x$ and~$y$. The same computations proving Lemma~\ref{ModComStrucFirst} show also here that~$\phi$ and~$\phi'$ commute and both have order~$2$. Therefore we get, as in Paragraph~\ref{ModComStrucFirst} and using the notation introduced there, an $H$-module algebra structure on~$A$ by requiring that~$g_2$ acts by~$\phi$ and~$g_3$ acts by~$\phi'$. 

From this action, we construct the coaction again in exactly the same way as in Paragraph~\ref{ModComStrucFirst}: We use the bicharacter~$\theta$ defined there to introduce the coaction via the formula
\[\delta \coln A \to H \ot A,~a \mapsto \frac{\K}{4} \sum_{g,g' \in G} \theta(g,g') \, g \ot g'.a\]
Because the formulas are the same, the same arguments as in Proposition~\ref{ModComStrucFirst} show that~$\ua$ and~$x^2$ are coinvariant and that 
\begin{align*}
&\delta(x) = \frac{\K}{2} (g_1 + g_3) \ot x + \frac{\K}{2} (g_1 - g_3) \ot x^3 \\
&\delta(x^3) = \frac{\K}{2} (g_1 - g_3) \ot x + \frac{\K}{2} (g_1 + g_3) \ot x^3
\end{align*}

\subsection[The coalgebra structure of the second example]{} \label{CoalgStrucSec}
As in Paragraph~\ref{CoalgStrucFirst}, we have now introduced sufficiently many structure elements to form the algebra~$A \hat{\ot} A$, and we define the coproduct as the unique algebra homomorphism 
$\da \coln A \to A \hat{\ot} A$ that maps~$x$ to~$x'$ and~$y$ to~$y'$, where
\[x' \deq  \frac{\K}{2} (x \ot x + x \ot x^3 + x^3 \ot x - x^3 \ot x^3)
\qquad \text{and} \qquad y' \deq y \ot y\]
To see that there is such an algebra homomorphism, we need to show that~$x'$ and~$y'$ satisfy the same relations as~$x$ and~$y$, i.e., the defining relations of~$A$:
\begin{prop}
\[x'^4 = 1 \qquad \qquad x'y' = y'x'^3 \qquad \qquad 
y'^2 = \frac{\K}{2}(\zeta 1 + x' - \zeta x'^2 + x'^3)\]
\end{prop}
\begin{pf}
\begin{pflist}
\item
Note that the elements~$x'$ and~$y'$ have the same form as the ones in Paragraph~\ref{CoalgStrucFirst}, although they lie in a different algebra. However, the algebra structure on the subalgebra spanned by the elements of the form~$x^i \ot x^j$, for 
$i,j=0,1,2,3$, is actually the same. Therefore, the equation
\[(\ua \ot x) (x \ot \ua) = x \ot x = (x \ot \ua) (\ua \ot x)\]
found in Lemma~\ref{CoalgStrucFirst} also holds in our situation, and as in the proof of Proposition~\ref{CoalgStrucFirst}, we have $x'^2 = x^2 \ot x^2$ and~$x'^4 = 1$.

\item
Again as in the proof of Proposition~\ref{CoalgStrucFirst}, we now compute the square of~$y'$ in~$A \hat{\ot} A$. In the present case, the relation
$y^2 = \frac{\K}{2}(x (\ua + x^2) + \zeta (\ua - x^2))$ implies
\[(\ua + x^2) y^2 = x + x^3 \qquad \text{and} \qquad 
(\ua - x^2) y^2 = \zeta (\ua - x^2)\]
Exactly as before, we have
\begin{align*}
&y'^2 = \frac{\K}{4} \sum_{i,j=1}^4 \theta(g_i,g_j) y (g_i.y) \ot (g_j.y) y 
\end{align*}
Again, we look at the terms with~$i$ fixed individually, and carry out the summation over~$j$. For $i=1$, we have
\begin{align*}
\sum_{j=1}^4 \theta(g_1,g_j) (g_j.y) y &= \sum_{j=1}^4 (g_j.y) y = 
y^2 + x^3 y^2 + x^2 y^2 + x y^2 = \ua + x + x^2 + x^3   
\end{align*}

For $i=2$, we have
\begin{align*}
&\sum_{j=1}^4 \theta(g_2,g_j) (g_j.y) y = 
y^2 + \zeta^2 x^3 y^2 - x^2 y^2 - \zeta^2 x y^2 \\
&= (\ua - x^2) y^2 - \zeta^2 x (\ua - x^2) y^2 
= \zeta (\ua - x^2) - \zeta^3 x (\ua - x^2) \\
&= \zeta \ua - \zeta^3 x - \zeta x^2 + \zeta^3 x^3
\end{align*}

For $i=3$, we have
\begin{align*}
&\sum_{j=1}^4 \theta(g_3,g_j) (g_j.y) y = 
y^2 - x^3 y^2 + x^2 y^2 - x y^2 = - \ua + x - x^2 + x^3
\end{align*}

For $i=4$, we have
\begin{align*}
&\sum_{j=1}^4 \theta(g_4,g_j) (g_j.y) y = 
y^2 - \zeta^2 x^3 y^2 - x^2 y^2 + \zeta^2 x y^2 \\
&= (\ua - x^2) y^2 + \zeta^2 x (\ua - x^2) y^2 
= \zeta (\ua - x^2) + \zeta^3 x (\ua - x^2) \\
&= \zeta \ua + \zeta^3 x - \zeta x^2 - \zeta^3 x^3 
\end{align*}

In total, we obtain
\begin{align*}
4 y'^2 &= y^2 \ot (\ua + x + x^2 + x^3) +
y x^3 y \ot (\zeta \ua - \zeta^3 x - \zeta x^2 + \zeta^3 x^3)  \\
&+ y x^2 y \ot (- \ua + x - x^2 + x^3)  +
y x y \ot (\zeta \ua + \zeta^3 x - \zeta x^2 - \zeta^3 x^3)  
\end{align*}
By the commutation relation between~$x$ and~$y$ in~$A$, this is equal to
\begin{align*}
4 y'^2 &= y^2 \ot (\ua + x + x^2 + x^3) +
x y^2 \ot (\zeta \ua - \zeta^3 x - \zeta x^2 + \zeta^3 x^3)  \\
&+ x^2 y^2 \ot (- \ua + x - x^2 + x^3)  +
x^3 y^2 \ot (\zeta \ua + \zeta^3 x - \zeta x^2 - \zeta^3 x^3)  
\end{align*}
Collecting terms, this becomes
\begin{align*}
4 y'^2 &= (\ua + \zeta x - x^2 + \zeta x^3) y^2 \ot \ua +
(\ua - \zeta^3 x + x^2 + \zeta^3 x^3) y^2 \ot x \\
&+ (\ua - \zeta x - x^2 - \zeta x^3) y^2 \ot x^2   +
(\ua + \zeta^3 x + x^2 - \zeta^3 x^3) y^2 \ot x^3 
\end{align*}
or, with the help of the equations recorded above,
\begin{align*}
4 y'^2 &= (\zeta \ua + \zeta \ua - \zeta x^2 + \zeta x^2)  \ot \ua +
(x - \zeta^4 x + x^3 + \zeta^4 x^3) \ot x \\
&+ (\zeta \ua - \zeta \ua  - \zeta x^2 - \zeta x^2) \ot x^2   +
(x + \zeta^4 x + x^3 - \zeta^4 x^3) \ot x^3 
\end{align*}
which simplifies to
\begin{align*}
y'^2 &= \frac{\K}{2} (\zeta \ua \ot \ua + x^3 \ot x - \zeta x^2  \ot x^2 + x \ot x^3) 
\end{align*}

On the other hand, we have
\begin{align*}
\frac{\K}{2}(&\zeta 1 + x' - \zeta x'^2 + x'^3) = \\
&\frac{\K}{2}(\zeta \ua \ot \ua + \frac{\K}{2} (x \ot x + x \ot x^3 + x^3 \ot x - x^3 \ot x^3) \\
&\mspace{20mu}
- \zeta x^2 \ot x^2 + \frac{\K}{2} (x^3 \ot x^3 + x^3 \ot x + x \ot x^3 - x \ot x)) \\
&= 
\frac{\K}{2}(\zeta \ua \ot \ua - \zeta x^2 \ot x^2 + x \ot x^3 + x^3 \ot x)
\end{align*}
which establishes that $y'^2 = \frac{\K}{2}(\zeta 1 + x' - \zeta x'^2 + x'^3)$.

\item
We claim that $(x \ot x)(y \ot y) = (y \ot y)(x^3 \ot x^3)$. As we said at the end of Paragraph~\ref{ModComStrucSec}, the formulas for the coaction from Proposition~\ref{ModComStrucFirst} still hold in our situation. Therefore, we obtain exactly as in the third step of the proof of Proposition~\ref{CoalgStrucFirst} that
\begin{align*}
(\ua \ot x)(y \ot \ua) &= \frac{\K}{2} (\ua + x^2) y \ot x + \frac{\K}{2} (\ua - x^2) y \ot x^3
\end{align*}
and
\begin{align*}
(\ua \ot y)(x \ot \ua) &= \frac{\K}{2} x \ot (\ua + x^2) y + \frac{\K}{2} x^3 \ot (\ua - x^2)y
\end{align*}
From the first of these two formulas, we get
\begin{align*}
(x \ot x)(y \ot y) 
&= \frac{\K}{2} (x + x^3) y \ot x y + \frac{\K}{2} (x - x^3) y \ot x^3 y \\
&= \frac{\K}{2} (xy  \ot xy + x^3y \ot xy + xy \ot x^3 y - x^3 y \ot x^3 y)
\end{align*}
while the second of these formulas implies that
\begin{align*}
(y \ot y)(x \ot x) &= 
\frac{\K}{2} yx \ot (\ua + x^2) y x + \frac{\K}{2} yx^3 \ot (\ua - x^2)yx \\
&= 
\frac{\K}{2} (yx \ot yx + yx \ot x^2 yx + yx^3 \ot y x - yx^3 \ot x^2 y x)
\end{align*}

The element~$x^2$ is invariant, coinvariant, and central in~$A$. Therefore, if we multiply the previous formula by $x^2 \ot x^2$, we obtain
\begin{align*}
(y \ot y)(x^3 \ot x^3) &= 
\frac{\K}{2} (yx^3 \ot yx^3 + yx^3 \ot yx + yx \ot y x^3 - yx \ot y x)
\end{align*}
Now the commutation relation between~$x$ and~$y$ yields our claim. 

\item
The properties of~$x^2$ just stated imply that the elements $x^2 \ot \ua$ and $\ua \ot x^2$ are central in~$A \hat{\ot} A$. The claim established in the previous step therefore also yields that
\[(x \ot x^3)(y \ot y) = (y \ot y)(x^3 \ot x) \quad \text{and} \quad 
(x^3 \ot x)(y \ot y) = (y \ot y)(x \ot x^3)\]
as well as $(x^3 \ot x^3)(y \ot y) = (y \ot y)(x \ot x)$. We therefore get
\begin{align*}
x'y' &=  \frac{\K}{2} (x \ot x + x \ot x^3 + x^3 \ot x - x^3 \ot x^3)(y \ot y) \\
&= \frac{\K}{2} (y \ot y)(x^3 \ot x^3 + x^3 \ot x + x \ot x^3 - x \ot x) \\
&= \frac{\K}{2} (y \ot y)(x^2 \ot x^2) (x \ot x + x \ot x^3 + x^3 \ot x - x^3 \ot x^3)
\end{align*}
Because we have already seen in the first step that~$x'^2 = x^2 \ot x^2$, this proves that
$x'y' = y'x'^3$, as asserted.
\qed
\end{pflist}
\end{pf}

The very same computations that we have carried out in Paragraph~\ref{CoalgStrucFirst} show that $(\phi \ot \phi) \circ \da = \da \circ \phi$ and 
$(\phi' \ot \phi') \circ \da = \da \circ \phi'$. As explained there, this implies that~$\da$ is $H$-linear and colinear.

We define the counit to be the unique algebra homomorphism~$\ea \coln A \to K$ which takes the values
\[\ea(x) = \K = \ea(y)\]
on the generators. Because the images of the generators satisfy the defining relations of our algebra, $\ea$ is well-defined. It is also immediate to check that
$\ea \circ \phi = \ea = \ea \circ \phi'$, which implies as in Paragraph~\ref{CoalgStrucFirst} that $\ea$ is $H$-linear and colinear. The same computations that we used there to show that
\[(\ea \ot \id_A)(\da(a)) = \K \ot a \qquad \text{and} \qquad 
(\id_A \ot \, \ea)(\da(a)) = a \ot \K\]
also establish these equations in the present case.

\subsection[A basis consisting of group-like elements for the second example]{} \label{BasisGrouplSec}
Again, we prove the coassociativity of the coproduct by exhibiting a basis consisting of group-like elements. They are given by the same formulas as in Paragraph~\ref{BasisGrouplFirst}:
\begin{defn}
We define $\omega_1 \deq \ua$,
\[\omega_2 \deq \frac{\K}{2}(1 + \iota \zeta^2) x + \frac{\K}{2}(1 - \iota \zeta^2) x^3 \qquad
\omega_3 \deq \frac{\K}{2}(1 - \iota \zeta^2) x + \frac{\K}{2}(1 + \iota \zeta^2) x^3\] 
and $\omega_4 \deq x^2$. Furthermore, we define
\[\eta_1 \deq y \qquad \eta_2 \deq x^3 y \qquad \eta_3 \deq x^2 y \qquad \eta_4 \deq x y\]
\end{defn}

As before, these elements have the required properties:
\begin{prop}
The elements $\omega_1, \omega_2, \omega_3, \omega_4, \eta_1, \eta_2, \eta_3, \eta_4$ form a basis of~$A$. They are group-like. 
\end{prop}
\begin{pf}
The same argument as in Paragraph~\ref{BasisGrouplFirst} shows that these elements form a basis. Because the coproduct in the first example and the coproduct in the present second example agree on the space~$\Span(\ua, x, x^2, x^3)$, it follows from our previous computation that the elements $\omega_1$, $\omega_2$, $\omega_3$, and~$\omega_4$ are group-like. 

The element~$\eta_1 \deq y$ is group-like by construction. Because~$\eta_2 \deq x^3 y = \phi(y)$, it was already necessary to show that $\eta_2$ is group-like when establishing the $H$-linearity of the coproduct. But then~$\eta_3 = \omega_4 \eta_1$ and~$\eta_4 = \omega_4 \eta_2$ are also group-like, because the product of a coinvariant group-like element with an arbitrary group-like element is again group-like, a fact already mentioned in the proof of Proposition~\ref{BasisGrouplFirst}.
\qed
\end{pf}

\begin{samepage}
The products~$\omega_i \omega_j$ are clearly the same as the corresponding ones in Paragraph~\ref{BasisGrouplFirst}:
\begin{center}
\renewcommand\arraystretch{1.5}
\begin{tabular}{|c|c|c|c|c|}
\hline 
& $\omega_1$ & $\omega_2$ & $\omega_3$ & $\omega_4$   \\ 
\hline 
$\omega_1$ & $\omega_1$ & $\omega_2$ & $\omega_3$ & $\omega_4$   \\ 
\hline 
$\omega_2$ & $\omega_2$ & $\omega_1$ & $\omega_4$ & $\omega_3$ \\ 
\hline 
$\omega_3$ & $\omega_3$ & $\omega_4$ & $\omega_1$ & $\omega_2$ \\ 
\hline 
$\omega_4$ & $\omega_4$ & $\omega_3$ & $\omega_2$ & $\omega_1$ \\ 
\hline 
\end{tabular} 
\end{center}
\end{samepage}

The products~$\omega_i \eta_j$ are also the same, because the commutation relation 
between~$x$ and~$y$ is not used in their computation:
\begin{center}
\renewcommand\arraystretch{1.5}
\begin{tabular}{|c|c|c|c|c|}
\hline 
& $\omega_1$ & $\omega_2$ & $\omega_3$ & $\omega_4$   \\ 
\hline 
$\eta_1$ & $\eta_1$ & 
$\frac{\K}{2}(1 - \iota \zeta^2) \eta_2 + \frac{\K}{2} (1 + \iota \zeta^2) \eta_4$ & $\frac{\K}{2}(1 + \iota \zeta^2) \eta_2 + \frac{\K}{2} (1 - \iota \zeta^2) \eta_4$ & 
$\eta_3$ \\ 
\hline 
$\eta_2$ & $\eta_2$ & 
$\frac{\K}{2}(1 + \iota \zeta^2) \eta_1 + \frac{\K}{2} (1 - \iota \zeta^2) \eta_3$ & 
$\frac{\K}{2}(1 - \iota \zeta^2) \eta_1 + \frac{\K}{2} (1 + \iota \zeta^2) \eta_3$ & $\eta_4$ \\ 
\hline 
$\eta_3$ & $\eta_3$ & 
$\frac{\K}{2}(1 + \iota \zeta^2) \eta_2 + \frac{\K}{2} (1 - \iota \zeta^2) \eta_4$ & $\frac{\K}{2}(1 - \iota \zeta^2) \eta_2 + \frac{\K}{2} (1 + \iota \zeta^2) \eta_4$ & $\eta_1$ \\ 
\hline 
$\eta_4$ & $\eta_4$ & 
$\frac{\K}{2}(1 - \iota \zeta^2) \eta_1 + \frac{\K}{2} (1 + \iota \zeta^2) \eta_3$ & 
$\frac{\K}{2}(1 + \iota \zeta^2) \eta_1 + \frac{\K}{2} (1 - \iota \zeta^2) \eta_3$ & $\eta_2$ \\ 
\hline 
\end{tabular} 
\end{center}

However, as our present algebra is noncommutative, we do not have~$\omega_i \eta_j = \eta_j \omega_i$ in general; this equation holds only if~$i=1$ or~$i=4$, because~$\omega_1$ and~$\omega_4$ are central. For~$i=2$ and~$i=3$, the commutation relation between~$x$ and~$y$ implies immediately that~$\omega_2 \eta_j = \eta_j \omega_3$ and~$\omega_3 \eta_j = \eta_j \omega_2$, so that the products~$\eta_j \omega_i$ are given by the following table:
\begin{center}
\renewcommand\arraystretch{1.5}
\begin{tabular}{|c|c|c|c|c|}
\hline 
& $\omega_1$ & $\omega_2$ & $\omega_3$ & $\omega_4$   \\ 
\hline 
$\eta_1$ & $\eta_1$ & 
$\frac{\K}{2}(1 + \iota\zeta^2) \eta_2 + \frac{\K}{2} (1 - \iota \zeta^2) \eta_4$ & $\frac{\K}{2}(1 - \iota\zeta^2) \eta_2 + \frac{\K}{2}(1+ \iota\zeta^2) \eta_4$& 
$\eta_3$ \\ 
\hline 
$\eta_2$ & $\eta_2$ & 
$\frac{\K}{2}(1 - \iota \zeta^2) \eta_1 + \frac{\K}{2} (1 + \iota \zeta^2) \eta_3$ & 
$\frac{\K}{2}(1 + \iota \zeta^2) \eta_1 + \frac{\K}{2} (1 - \iota \zeta^2) \eta_3$ & $\eta_4$ \\ 
\hline 
$\eta_3$ & $\eta_3$ & 
$\frac{\K}{2}(1 - \iota \zeta^2) \eta_2 + \frac{\K}{2} (1 + \iota \zeta^2) \eta_4$ & 
$\frac{\K}{2}(1 + \iota \zeta^2) \eta_2 + \frac{\K}{2} (1 - \iota \zeta^2) \eta_4$ & $\eta_1$ \\ 
\hline 
$\eta_4$ & $\eta_4$ & 
$\frac{\K}{2}(1 + \iota \zeta^2) \eta_1 + \frac{\K}{2} (1 - \iota \zeta^2) \eta_3$ & 
$\frac{\K}{2}(1 - \iota \zeta^2) \eta_1 + \frac{\K}{2} (1 + \iota \zeta^2) \eta_3$ & $\eta_2$ \\ 
\hline 
\end{tabular} 
\end{center}

The products $\eta_i \eta_j$ can be computed relatively easily from the fact that
\begin{align*}
\eta_1 \eta_1 = y^2 = \frac{\K}{2}(\zeta \ua + x - \zeta x^2 + x^3) =
\frac{\K}{2}(\zeta \omega_1 + \omega_2 + \omega_3 -\zeta \omega_4)
\end{align*}
together with the equation 
$- \frac{\iota}{\zeta} \omega_2 + \frac{\iota}{\zeta} \omega_3 = \zeta x - \zeta x^3$ encountered already in Paragraph~\ref{BasisGrouplFirst}. They are given in the following table:
\begin{center}
\renewcommand\arraystretch{1.5}
\begin{tabular}{|c|c|c|}
\hline  
 & $\eta_1$ & $\eta_2$  \\ 
\hline 
$\eta_1$ & 
$\frac{\zeta}{2} \omega_1 + \frac{\K}{2} \omega_2 
+ \frac{\K}{2} \omega_3 - \frac{\zeta}{2} \omega_4$ & 
$\frac{\K}{2} \omega_1 - \frac{\iota}{2\zeta} \omega_2 
+ \frac{\iota}{2\zeta} \omega_3 + \frac{\K}{2} \omega_4$  \\ 
\hline 
$\eta_2$ & 
$\frac{\K}{2} \omega_1 + \frac{\iota}{2\zeta} \omega_2 
- \frac{\iota}{2\zeta} \omega_3 + \frac{\K}{2} \omega_4$ & 
$\frac{\zeta}{2} \omega_1 + \frac{\K}{2} \omega_2 
+ \frac{\K}{2} \omega_3 - \frac{\zeta}{2} \omega_4$  \\ 
\hline 
$\eta_3$ & 
$- \frac{\zeta}{2} \omega_1 + \frac{\K}{2} \omega_2 
+ \frac{\K}{2} \omega_3 + \frac{\zeta}{2} \omega_4$ & 
$\frac{\K}{2} \omega_1 + \frac{\iota}{2\zeta} \omega_2 
- \frac{\iota}{2\zeta} \omega_3 + \frac{\K}{2} \omega_4$  \\ 
\hline 
$\eta_4$ & 
$\frac{\K}{2} \omega_1 - \frac{\iota}{2\zeta} \omega_2 
+ \frac{\iota}{2\zeta} \omega_3 + \frac{\K}{2} \omega_4$ & 
$- \frac{\zeta}{2} \omega_1 + \frac{\K}{2} \omega_2 
+ \frac{\K}{2} \omega_3 + \frac{\zeta}{2} \omega_4$  \\ 
\hline 
\end{tabular} 
\end{center}

\begin{center}
\renewcommand\arraystretch{1.5}
\begin{tabular}{|c|c|c|}
\hline  
 & $\eta_3$ & $\eta_4$ \\ 
\hline 
$\eta_1$ & 
$- \frac{\zeta}{2} \omega_1 + \frac{\K}{2} \omega_2 
+ \frac{\K}{2} \omega_3 + \frac{\zeta}{2} \omega_4$ & 
$\frac{\K}{2} \omega_1 + \frac{\iota}{2\zeta} \omega_2 
- \frac{\iota}{2\zeta} \omega_3 + \frac{\K}{2} \omega_4$ \\ 
\hline 
$\eta_2$ & 
$\frac{\K}{2} \omega_1 - \frac{\iota}{2\zeta} \omega_2 
+ \frac{\iota}{2\zeta} \omega_3 + \frac{\K}{2} \omega_4$ & 
$- \frac{\zeta}{2} \omega_1 + \frac{\K}{2} \omega_2 
+ \frac{\K}{2} \omega_3 + \frac{\zeta}{2} \omega_4$ \\ 
\hline 
$\eta_3$ & 
$\frac{\zeta}{2} \omega_1 + \frac{\K}{2} \omega_2 
+ \frac{\K}{2} \omega_3 - \frac{\zeta}{2} \omega_4$  & 
$\frac{\K}{2} \omega_1 - \frac{\iota}{2\zeta} \omega_2 
+ \frac{\iota}{2\zeta} \omega_3 + \frac{\K}{2} \omega_4$ \\ 
\hline 
$\eta_4$ & 
$\frac{\K}{2} \omega_1 + \frac{\iota}{2\zeta} \omega_2 
- \frac{\iota}{2\zeta} \omega_3 + \frac{\K}{2} \omega_4$ & 
$\frac{\zeta}{2} \omega_1 + \frac{\K}{2} \omega_2 
+ \frac{\K}{2} \omega_3 - \frac{\zeta}{2} \omega_4$ \\ 
\hline 
\end{tabular} 
\end{center}
In this table, it is understood that the first factor~$\eta_i$ is listed in the first column, whereas the second factor~$\eta_j$ is listed in the first row. Note that, in contrast to the first example, this table is not symmetric, because we do not have $\eta_i \eta_j = \eta_j \eta_i$ in general. 

Because the action of the group~$G$ was introduced by the same formulas in both examples, its action on the group-like elements of~$A$ has the same form as in the first example, which was described in Lemma~\ref{BasisGrouplFirst}:
\begin{lem}
\leavevmode
\begin{thmlist}
\item 
$\omega_1$ and~$\omega_4$ are fixed points of the $G$-action. Moreover, we have $g_3.\omega_i = \omega_i$ for all $i=1,2,3,4$, while $g_2.\omega_2 = \omega_3$.

\item
$G$ acts simply transitively on the set~$\{\eta_1, \eta_2, \eta_3, \eta_4\}$: If we define
$\eta_{g_i} \deq \eta_i$, we have $g_i.\eta_{g_j} = \eta_{g_i g_j}$.
\end{thmlist}
\end{lem}

\subsection[The antipode in the second example]{} \label{AntipSec}
Also in our present case, the last structure element that we need to construct is the antipode. Here, we define it to be the unique $K$-linear endomorphism~$\sa \coln A \to A$ that is given on the basis consisting of group-like elements by the formulas 
\begin{align*}
\sa(\omega_1) &= \omega_1 \qquad \sa(\omega_4) = \omega_4 \quad &
\sa(\omega_2) &= \omega_2 \qquad \sa(\omega_3) = \omega_3 \\
\sa(\eta_1) &= 
\frac{\K}{2}\left(\frac{\K}{\zeta} \eta_1 + \eta_2 - \frac{\K}{\zeta} \eta_3 + \eta_4 \right) \quad &
\sa(\eta_2) &= 
\frac{\K}{2}\left(\eta_1 + \frac{\K}{\zeta} \eta_2 + \eta_3 - \frac{\K}{\zeta} \eta_4 \right)  \\
\sa(\eta_3) &= 
\frac{\K}{2}\left(- \frac{\K}{\zeta}\eta_1 + \eta_2 + \frac{\K}{\zeta} \eta_3 + \eta_4 \right) \quad &
\sa(\eta_4) &= 
\frac{\K}{2}\left(\eta_1 - \frac{\K}{\zeta} \eta_2 + \eta_3 + \frac{\K}{\zeta} \eta_4 \right)
\end{align*}

This endomorphism has the required properties:
\begin{prop}
$\sa$ is an antipode for~$A$.
\end{prop}
\begin{pf}
We have already explained in the proof of Proposition~\ref{AntipFirst} that it suffices to show that~$\sa$ maps each of our group-like basis elements to its inverse, and for the elements $\omega_1$, $\omega_2$, $\omega_3$, and~$\omega_4$, this holds for the same reasons already discussed there. For the basis element~$\eta_1$, we get from the multiplication table given in Paragraph~\ref{BasisGrouplSec} that
\begin{align*}
\eta_1 \sa(\eta_1) &= 
\frac{\K}{2} \Bigg(\frac{\K}{\zeta} \eta_1 \eta_1 + \eta_1 \eta_2 - \frac{\K}{\zeta} \eta_1 \eta_3 
+ \eta_1 \eta_4 \Bigg) \\
&= \frac{\K}{4} \Bigg(\frac{\K}{\zeta} (\zeta \omega_1 +  \omega_2 +  \omega_3 - \zeta \omega_4) 
+ (\omega_1 - \frac{\iota}{\zeta} \omega_2 
+ \frac{\iota}{\zeta} \omega_3 + \omega_4) \\
&- \frac{\K}{\zeta} (- \zeta \omega_1 + \omega_2 +  \omega_3 + \zeta \omega_4)
+  (\omega_1 + \frac{\iota}{\zeta} \omega_2 - \frac{\iota}{\zeta} \omega_3 + \omega_4)\Bigg)
= \omega_1
\end{align*}
As in the proof of Proposition~\ref{AntipFirst}, it would be possible to treat the remaining basis elements~$\eta_2$, $\eta_3$, and~$\eta_4$ in the same way, and then to use Lemma~\ref{YetDrinfHopf} to conclude that~$\sa$ is $H$-linear and colinear. However, as it was there, it is easier here to verify explicitly that~$\sa$ is $H$-linear, and therefore also colinear, and then to deduce the equation 
$\eta_i \sa(\eta_i) = \omega_1$ for $i=2,3,4$ by acting with~$g_i$ on the equation 
$\eta_1 \sa(\eta_1) = \omega_1$.
\qed
\end{pf}

\subsection[The algebra structure]{} \label{AlgStruc}
We have seen already in Proposition~\ref{AlgStrucFirst} that the algebra constructed for the first example is semisimple and already splits over the base field used. This also holds for the algebra~$A$ in our present, second example:
\begin{prop}
As algebras, we have~$A \cong K^4 \oplus M(2 \times 2, K)$.
\end{prop}
\begin{pf}
\begin{pflist}
\item 
In the proof of Proposition~\ref{AntipSec}, we saw explicitly that~$\eta_1 = y$ is invertible, so that a one-dimensional character must map it to a nonzero element. The relation $xy = yx^3$ then implies that~$x$ cannot be mapped to an arbitrary fourth root of unity, but only to~$\pm \K$. Depending on this value for~$x$, we have the following values for $\frac{\K}{2}(\zeta \ua + x - \zeta x^2 + x^3)$:
\begin{center}
\renewcommand\arraystretch{1.5}
\begin{tabular}{|c|c|c|} \hline
$x$ & $\K$ & $-\K$  \\ \hline
$\frac{\K}{2}(\zeta \ua + x - \zeta x^2 + x^3)$ & $\K$ & $-\K$  \\
\hline
\end{tabular}
\end{center}

\enlargethispage{6pt}
From this table, we see that there are indeed four one-dimensional representations $\varepsilon_1$, $\varepsilon_2$, $\varepsilon_3$, and~$\varepsilon_4$ of~$A$ that take the following values on the generators~$x$ and~$y$:
\begin{samepage}
\begin{center}
\renewcommand\arraystretch{1.5}
\begin{tabular}{|c|c|c|c|c|} \hline
& $\varepsilon_1$ & $\varepsilon_2$ & $\varepsilon_3$ & $\varepsilon_4$  \\ \hline
$x$ & $\K$ & $-\K$ & $-\K$ & $\K$ \\
\hline
$y$ & $\K$ & $\iota$ & $-\iota$ & $-\K$  \\
\hline
\end{tabular}
\end{center}
Note that we have $\varepsilon_1 = \ea$, the counit introduced in Paragraph~\ref{CoalgStrucSec}.
\end{samepage}

\item
Next, we exhibit an irreducible two-dimensional representation of~$A$. Proposition~\ref{AlgSec} implies that the elements $\ua, x, x^2, x^3, y, yx, yx^2, y x^3$ form a basis of~$A$; in fact, this basis is just a permutation of the basis given there. This shows that~$A$ is free as a right module over the subalgebra \mbox{$C \deq \Span(\ua, x, x^2, x^3)$}; the elements~$\ua$ and~$y$ form a basis for this module. By assigning the value~$\iota$ to~$x$, we get a one-dimensional representation not of~$A$, but of~$C$. As with all our one-dimensional representations, the underlying vector space is the base field~$K$. Therefore we can construct the induced $A$-module $V \deq A \ot_C K$ of this one-dimensional $C$-module. As a consequence of the freeness of~$A$ over~$C$ just discussed, the elements~$\ua \ot_C \K$ and~$y \ot_C \K$ form a basis for~$V$. It follows directly from the defining relations that the actions of the generators~$x$ and~$y$ with respect to this basis are represented by the $2 \times 2$-matrices
\[X = \begin{pmatrix} \iota & 0 \\ 0 & -\iota \end{pmatrix} \qquad \text{and} \qquad
Y = \begin{pmatrix} 0 & \zeta \\ 1 & 0 \end{pmatrix}\]
which should not be confused with the $8 \times 8$-matrices of the same name used in Paragraph~\ref{AlgSec}. From the form of these matrices, we can infer that~$V$ is an irreducible $A$-module, as a nonzero proper submodule would be one-dimensional and therefore be spanned by a vector that must correspond to an eigenvector for both~$X$ and~$Y$. Now an eigenvector for~$X$ must be a multiple of one of the unit vectors, but none of these multiples is also an eigenvector for~$Y$.

\item
Now consider the Jacobson radical~$J$ of~$A$. $A/J$ is a semisimple algebra of dimension at most~$8$ (cf.~\cite{FD}, Cor.~2.3, p.~60). The Jacobson radical is defined as the intersection of the annihilators of all the simple modules; in particular, it annihilates the four one-dimensional representations and the one two-dimensional representation constructed above, so that these are also modules for~$A/J$. The Wedderburn decomposition of~$A/J$ therefore contains four two-sided ideals that are isomorphic to~$K$ and one two-sided ideal that is isomorphic to~$M(2 \times 2, K)$. In particular, its dimension is at least~$8$. This shows that
\[\dm_K(A/J) = 8 = \dm_K(A)\] so~$J=\{0\}$ and~$A$ is semisimple. It also must have the asserted Wedderburn decomposition. 
\qed
\end{pflist}
\end{pf}

\section{Cores} \label{Sec:Cores}
\subsection[Products of group-like elements]{} \label{ProdGroupl}
In this section, we assume that our base field~$K$ is algebraically closed of characteristic zero. We assume that~$G$ is a finite abelian group, and we denote its group ring by~$H \deq K[G]$. The basic object of study will be a finite-dimensional cocommutative cosemisimple Yetter-Drinfel'd Hopf algebra~$A$ over~$H$. It follows from these assumptions that~$A$ has a basis consisting of group-like elements. In this paragraph and the next, we will review some results from~\cite{SoTriv} in a dualized setting, in order to discuss what they mean for our examples constructed in the previous sections. 

The action of~$G$ on~$A$ can be equivalently described by a representation that we denote by
\[G \to \Aut(A),~g \mapsto \phi_g\]
For the character group~$\hat{G} = \Hom(G,K^\times)$ of~$G$, there is a similar representation. Each element of the character group can be linearly extended to an algebra homomorphism from~$H$ to~$K$, and in this way~$\hat{G}$ can be regarded as a basis of~$H^*$, so that~$H^*$ is isomorphic to the group ring of~$\hat{G}$. The left coaction of~$H$ on~$A$ gives rise to a right action of~$H^*$ on~$A$, but as~$H^*$ is commutative, this can be viewed as a left module structure and therefore leads to a representation 
\[\hat{G} \to \Aut(A),~\gamma \mapsto \psi_\gamma\]
of the character group that is given explicitly by $\psi_\gamma(a) = \gamma(a^\1) a^\2$. Because, as mentioned in Paragraph~\ref{YetMod}, the Yetter-Drinfel'd condition reduces to the dimodule condition in our situation, $\phi_g$ and~$\psi_\gamma$ commute.

In analogy with~\cite{SoTriv}, Def.~3.1, p.~494, we now associate various groups and parameters with a given group-like element:
\begin{defn}
Suppose that~$\eta \in A$ is a group-like element. 
\begin{thmlist}
\item 
The inertia group of~$\eta$ is the group $T_\eta \deq \{g \in G \mid \phi_g(\eta) = \eta\}$. 

\item 
The isotropy group of~$\eta$ is the group $Q_\eta \deq \{\gamma \in \hat{G} \mid \psi_\gamma(\eta) = \eta\}$. 

\item 
The index group of~$\eta$ is the group $G_\eta \deq Q_\eta^\perp/(T_\eta \cap Q_\eta^\perp)$. 

\item 
The index of~$\eta$ is the number $|G_\eta|$. 

\item 
The orbit of~$\eta$ under the action of~$Q_\eta^\perp$ is the set $O_\eta \deq \{\phi_g(\eta) \mid g \in Q_\eta^\perp\}$.
\end{thmlist}
\end{defn}

We note that it follows from~\cite{SoTriv}, Lem.~2.2, p.~492 that the index can also be written in the form $|G_\eta| = |T_\eta^\perp/(Q_\eta \cap T_\eta^\perp)|$, so that the definition is symmetric with respect to~$G$ and~$\hat{G}$.

The preceding definition applies only to one group-like element. We now discuss what happens when two group-like elements interact. As we have stated already in the proof of Proposition~\ref{BasisGrouplFirst}, the product of two group-like elements is in general not again group-like, and the multiplication tables for the group-like elements in Paragraph~\ref{BasisGrouplFirst} and Paragraph~\ref{BasisGrouplSec} show this explicitly. To see what happens instead, we consider two group-like elements~$\eta$ and~$\eta'$ and define the orbits
\[O \deq \{\psi_\gamma(\eta) \mid \gamma \in T_{\eta'}^\perp\} \qquad \qquad
O' \deq \{\phi_g(\eta') \mid g \in Q_{\eta}^\perp\}\]
of~$T_{\eta'}^\perp$ and~$Q_{\eta}^\perp$, respectively. We can then say the following:
\begin{prop}
\leavevmode
\begin{thmlist}
\item
$O$ and $O'$ have the same cardinality~$m$.

\item
There are distinct group-like elements $\omega_1,\ldots,\omega_m$ such that 
$$\tilde{\eta} \tilde{\eta}' \in \Span(\omega_1,\ldots,\omega_m)$$
for all $\tilde{\eta} \in O$ and all $\tilde{\eta}' \in O'$. 
In particular,
$\eta \eta' \in \Span(\omega_1,\ldots,\omega_m)$.

\item
If we write $\eta \eta'$ as a linear combination of $\omega_1,\ldots,\omega_m$, all coefficients are nonzero. Therefore, $\Span(\omega_1,\ldots,\omega_m)$ is the smallest subspace that contains~$\eta \eta'$ and is spanned by group-like elements. 
\end{thmlist}
\end{prop}
\begin{pf}
\begin{pflist}
\item
The dual space~$A^*$ of~$A$ is a right-right Yetter-Drinfel'd Hopf algebra over~$H$ (cf.~\cite{SoDevEnvAlg}, Subsec.~2.6, p.~37; the corresponding part of~\cite{SoDefUnivEinh} contains further details). The right $H$-structures can be dualized themselves, so that we get a left-left Yetter-Drinfel'd Hopf algebra structure on~$A^*$ over~$H^* \cong K[\hat{G}]$ (cf.~\cite{SoTriv}, Par.~1.4, p.~482 and the references given there). $A^*$~is commutative; the dual basis of the basis of group-like elements of~$A$ is a basis of~$A^*$ consisting of orthogonal primitive idempotents. The character corresponding to a primitive idempotent of~$A^*$ is given by evaluation at the corresponding group-like element;
in other words, the (one-dimensional) characters of~$A^*$ correspond under the isomorphism
$A \cong A^{**}$ exactly to the group-like elements of~$A$. 

\item
We can therefore apply the results from~\cite{SoTriv}, Sec.~3. For the Yetter-Drinfel'd Hopf algebra denoted by~$A$ in~\cite{SoTriv}, we use the dual~$A^*$ of our algebra, and for the finite abelian group~$G$ in~\cite{SoTriv}, we use our character group~$\hat{G}$. Note that substituting~$\hat{G}$ for~$G$ interchanges the roles of~$T_\eta$ and~$Q_\eta$. Our first assertion now follows from~\cite{SoTriv}, Lem.~3.2, p.~497, the second assertion from~\cite{SoTriv}, Thm.~3.3, p.~501, and the third assertion from the remark after that theorem.
\qed
\end{pflist}
\end{pf}

Because the issue was at least touched in Paragraph~\ref{BasisGrouplFirst} and Paragraph~\ref{BasisGrouplSec}, we also record what~\cite{SoTriv}, Cor.~3.3, p.~503 yields in our situation:
\begin{cor}
For two group-like elements~$\eta$ and~$\eta'$ in~$A$, the following conditions are equivalent:
\begin{enumerate}
\item 
$\eta \eta'$ is group-like.

\item
$Q_\eta^\perp \subset T_{\eta'}$

\item
$\sigma_{A,A}(\eta \ot \eta') = \eta' \ot \eta$
\end{enumerate}
\end{cor}

\subsection[The core of a group-like element]{} \label{CoreGroupl}
The definition of the core of a character can be found in~\cite{SoTriv}, Def.~3.8, p.~515. Under the dualization that we have considered in the preceding paragraph, this becomes the definition of the core of a group-like element, which we review now. 

So, suppose that~$\eta \in A$ is a group-like element. By~\cite{SoTriv}, Lem.~3.5, p.~508, $\sa(\eta)$ is in general not again group-like; this happens if and only if the index of~$\eta$ is~$1$. But in any case, because the group-like elements form a basis of~$A$, we can choose a group-like element~$\eta'$ that appears in~$\sa(\eta)$ with nonzero coefficient. If we apply Proposition~\ref{ProdGroupl} to this pair of group-like elements, the arising new group-like elements $\omega_1,\ldots,\omega_m$ have additional properties:
\begin{samepage}
\begin{prop}
\leavevmode
\begin{thmlist}
\item
We have $T_\eta = T_{\eta'}$ and $Q_\eta = Q_{\eta'}$.

\item
For all $i=1,\ldots,m$, we have $T_\eta \subset T_{\omega_i}$ and $Q_\eta \subset Q_{\omega_i}$.

\item 
For $\gamma \in T_\eta^\perp$, $\{\omega_1,\ldots,\omega_m\}$ is stable under~$\psi_\gamma$.

\item
For $g \in Q_\eta^\perp$, $\{\omega_1,\ldots,\omega_m\}$ is stable under~$\phi_g$.

\item
For some $i \le m$, we have $\omega_i = \ua$.

\item
$\Span(\omega_1,\ldots,\omega_m)$ is a subalgebra of~$A$.

\item
$\Span(\omega_1,\ldots,\omega_m)$ is stable under the antipode~$\sa$.

\item
$\{\phi_g(\eta) | g \in Q_\eta^\perp\} = \{\psi_\gamma(\eta) | \gamma \in T_\eta^\perp\}$

\item
$\Span(\{\phi_g(\eta) | g \in Q_\eta^\perp\}) = \Span(\omega_1 \eta,\ldots,\omega_m \eta)
= \Span(\eta \omega_1,\ldots,\eta \omega_m)$
\end{thmlist}
\end{prop}
\end{samepage}
\begin{pf}
We apply the results from~\cite{SoTriv}, Par.~3.7 and~3.8 to~$A^*$, considered as a left-left Yetter-Drinfel'd Hopf algebra over~$H^* \cong K[\hat{G}]$, as in the proof of Proposition~\ref{ProdGroupl}. We denote the dual basis elements of~$\eta$ and~$\eta'$ by~$e$ and~$e'$, respectively. Then~$e$ and~$e'$ are primitive idempotents in~$A^*$. In view of~\cite{SoTriv}, Cor.~3.5, p.~507, it follows from our choice of~$\eta'$ that~$e'$ is contained in the ideal~$\sa^*(I_e)$ introduced in~\cite{SoTriv}, Par.~3.1, p.~494, so that the assumptions made in~\cite{SoTriv}, Par.~3.7, p.~511 are satisfied. The first assertion then follows from the observations made at the beginning of that paragraph. Our second, third, fourth, and fifth assertion follow from the first, second, third, and fourth assertion in~\cite{SoTriv}, Prop.~3.7, p.~511. Our sixth assertion follows from~\cite{SoTriv}, Thm.~3.8, p.~515, our seventh assertion from~\cite{SoTriv}, Lem.~3.7, p.~512, our eighth assertion from~\cite{SoTriv}, Prop.~3.5, p.~507, and our ninth assertion from~\cite{SoTriv}, Thm.~3.7, p.~513.
\qed
\end{pf}

It also follows from the discussion in~\cite{SoTriv}, Par.~3.8, p.~514 that the group-like elements~$\omega_1,\ldots,\omega_m$ do not depend on the specific element~$\eta'$ that we choose, as long as its coefficient in the expansion of~$\sa(\eta)$ in the basis consisting of the group-like elements is nonzero. As a consequence, $\omega_1,\ldots,\omega_m$ depend on~$\eta$ alone. The following definition, which is the dualization of~\cite{SoTriv}, Def.~3.8, p.~515, is therefore meaningful:
\begin{defn}
For a group-like element~$\eta \in A$, the subalgebra $\Span(\omega_1,\ldots,\omega_m)$   of~$A$ is called the core of~$\eta$.
\end{defn}

It should be noted that the core is in general not a Yetter-Drinfel'd Hopf subalgebra of~$A$, because it is in general not an $H$-submodule; it is only stable under~$\phi_g$ 
for~$g \in Q_\eta^\perp$, not for an arbitrary~$g \in G$. It can be shown, however, that the core can be considered as a left-left Yetter-Drinfel'd Hopf algebra over the group ring~$K[G_\eta]$ of the index group (cf.~\cite{SoTriv}, Thm.~3.8, p.~515).

\subsection[The core in our examples]{} \label{CoreExamp}
We now consider the examples constructed in Section~\ref{Sec:FirstExamp} and Section~\ref{Sec:SecExamp} in the light of the theory recalled in a dualized form in the two preceding paragraphs. Before we can do that, we need to understand the relation between the mappings~$\phi_g$ and~$\psi_\gamma$ if the coaction is given by a nondegenerate bicharacter, as discussed in Paragraph~\ref{Bichar}:
\begin{lem}
For $g \in G$, define $\gamma \in \hat{G}$ by~$\gamma(h) \deq \theta(g,h)$ for all~$h \in G$.
Then we have~$\psi_\gamma = \phi_{g^{-1}}$.
\end{lem}
\begin{pf}
This follows from the form of the coaction described in Paragraph~\ref{QuasitriHopf} and Paragraph~\ref{Bichar}: We have
\begin{align*}
\psi_\gamma(a) &= \gamma(a^\1) a^\2 = 
\frac{\K}{|G|} \sum_{g',h' \in G} \theta(g',h') \, \gamma(h') \, g'.a \\
&= \frac{\K}{|G|} \sum_{g',h' \in G} \theta(g',h') \, \theta(g,h') \, g'.a
= \frac{\K}{|G|} \sum_{g',h' \in G} \theta(gg',h') \, g'.a
\end{align*}
Using the orthogonality relations from Lemma~\ref{Bichar}, we can carry out first the summation over~$h'$ and then over~$g'$ to obtain the assertion. 
\qed
\end{pf}

We discuss the first and the second example simultaneously. For the group-like element~$\eta$ whose core we considered in Paragraph~\ref{CoreGroupl}, we choose the 
element~\mbox{$\eta \deq \eta_1$}. From Lemma~\ref{BasisGrouplFirst} and Lemma~\ref{BasisGrouplSec}, we then see that~$T_\eta$ consists only of the unit element. By the lemma above, this implies that also~$Q_\eta$ consists only of the unit element, so that~$Q_\eta^\perp = G$ 
and~$T_\eta^\perp = \hat{G}$. Therefore, $\eta$ has index~$m=4$, and its core must be four-dimensional. The orbit of~$\eta$ in the sense of Definition~\ref{ProdGroupl} is the set
$\{\eta_1, \eta_2, \eta_3, \eta_4\}$. 

According to Paragraph~\ref{CoreGroupl}, we need to choose a group-like element~$\eta'$ that appears in~$\sa(\eta)$ with a nonzero coefficient. The formulas for the antipode in Paragraph~\ref{AntipFirst} and Paragraph~\ref{AntipSec} show that in both cases we can choose~$\eta' \deq \eta_1$. The formulas for the product~$\eta_1 \eta_1$ in Paragraph~\ref{BasisGrouplFirst} and Paragraph~\ref{BasisGrouplSec} then show that the core is spanned by the group-like elements $\omega_1$, $\omega_2$, $\omega_3$, and~$\omega_4$, as suggested by our notation. Moreover, the multiplication tables show that all products~$\eta_i \eta_j$ are contained in the core, in accordance with Proposition~\ref{ProdGroupl}. The multiplication tables also show that
\[\Span(\omega_1 \eta, \omega_2 \eta, \omega_3 \eta, \omega_4 \eta) = 
\Span(\eta \omega_1, \eta \omega_2, \eta \omega_3, \eta \omega_4) = 
\Span(\eta_1, \eta_2, \eta_3, \eta_4)\]
in accordance with Proposition~\ref{CoreGroupl}.

The multiplication table for the products~$\omega_i \omega_j$, which is the same in both examples, shows that the core is trivial in the sense of Definition~\ref{YetDrinfHopf}, because it is isomorphic to the group algebra~$K[G] = K[\Z_2 \times \Z_2]$ and therefore an ordinary Hopf algebra. However, in contrast to what happens in the prime order case (cf.~\cite{SoDiss}, Prop.~5.7, p.~68; \cite{SoYp}, Prop.~6.7, p.~98), it is not completely trivial: As we saw in Lemma~\ref{BasisGrouplFirst} and Lemma~\ref{BasisGrouplSec}, we have~$g_2.\omega_2 = \omega_3$. Moreover, it follows easily from Proposition~\ref{ModComStrucFirst} and the analogous formulas at the end of Paragraph~\ref{ModComStrucSec} that
\begin{align*}
\delta(\omega_2) 
&= \frac{\K}{2} (g_1 + g_3) \ot \omega_2 + \frac{\K}{2} (g_1 - g_3) \ot \omega_3 \\
\delta(\omega_3) 
&= \frac{\K}{2} (g_1 - g_3) \ot \omega_2 + \frac{\K}{2} (g_1 + g_3) \ot \omega_3
\end{align*}
in both examples.

Related to the fact that the core is not completely trivial, it is in our examples not the case
that the sets 
\[\{\omega_1 \eta, \omega_2 \eta, \omega_3 \eta, \omega_4 \eta\} \qquad \text{and} \qquad
\{\eta \omega_1, \eta \omega_2, \eta \omega_3, \eta \omega_4\} \]
are equal to~$\{\eta_1, \eta_2, \eta_3, \eta_4\}$; only the subspaces that they span coincide. This constitutes another difference to the prime order case cited above. In our examples, $\omega_i \eta$ is not always again group-like, as it was in the prime order case; rather, we see from the multiplication tables that for example
\[\omega_2 \eta_1 = 
\frac{\K}{2}(1 - \iota\zeta^2) \eta_2 + \frac{\K}{2} (1 + \iota \zeta^2) \eta_4\]
is a linear combination of two distinct group-like elements. This is in accordance with Proposition~\ref{ProdGroupl}: As we already stated above, we have~$T_\eta^\perp = \hat{G}$, and therefore find for the number~$m$ that appears in this proposition that
\[m = |\{\psi_\gamma(\omega_2) \mid \gamma \in T_\eta^\perp\}| = |\{\omega_2, \omega_3\}| = 2\]
where we have used the lemma above.

Although this is not quite as interesting, we can also compute the cores of~$\omega_2$ and~$\omega_3$. As there is no essential difference, we will focus on~$\omega_2$. From Lemma~\ref{BasisGrouplFirst} and Lemma~\ref{BasisGrouplSec}, we see that
$T_{\omega_2} = \{g_1,g_3\}$. As we have~$\theta(g_3,g_3) = 1$, but $\theta(g_2,g_3) = -1$, we can use the lemma above to obtain that also $Q_{\omega_2}^\perp = \{g_1,g_3\}$. Therefore, the index of~$\omega_2$ is $|Q_{\omega_2}^\perp/(T_{\omega_2} \cap Q_{\omega_2}^\perp)| = 1$. Its core is therefore one-dimensional. As the unit element~$\omega_1 = \ua$ is always contained in the core by Proposition~\ref{CoreGroupl}, we see that the core of~$\omega_2$ consists only of the scalar multiples of the unit element. Clearly, this core is completely trivial.

\subsection[The core-triviality conjecture]{} \label{CoreTrivConj}
In the examples that we have considered in this article, we have seen that the cores were trivial in the sense of Definition~\ref{YetDrinfHopf}. Although these examples certainly constitute only very little evidence, we conjecture that we have encountered here a general phenomenon:
\begin{conj}
Suppose that the base field~$K$ is algebraically closed of characteristic zero, and that~$A$ is a finite-dimensional cocommutative cosemisimple Yetter-Drinfel'd Hopf algebra over the group 
ring~$H \deq K[G]$ of a finite abelian group~$G$. Then the core of every group-like element~$\eta \in A$ is trivial as a Yetter-Drinfel'd Hopf algebra over the group ring of the index group~$G_\eta$.
\end{conj}
Additional, so far unpublished work of the authors provides further evidence for this conjecture.

\end{document}